\documentclass[a4paper,oneside,11pt,reqno]{amsart}

\usepackage{verbatim,upref,amsxtra,amssymb,amscd}

\usepackage{varioref}

\ifx\pdfoutput\undefined \usepackage[pagebackref,hypertex]{hyperref} \else \usepackage[pagebackref,pdftex]{hyperref} \fi

 \newcommand{\RR}{\mathbb{R}}        \newcommand{\NN}{\mathbb{N}}   \newcommand{\ZZ}{\mathbb{Z}}   

\DeclareMathAlphabet\mathscr{U}{eus}{m}{n} \SetMathAlphabet\mathscr{bold}{U}{eus}{b}{n} \DeclareMathAlphabet\matheur{U}{eur}{m}{n} \SetMathAlphabet\matheur{bold}{U}{eur}{b}{n}

\numberwithin{equation}{section}

\newtheorem{theo}{Theorem}[section] \newtheorem{prop}[theo]{Proposition} \newtheorem{lemm}[theo]{Lemma} \newtheorem{coro}[theo]{Corollary} 
\newtheorem{ques}{Question}

\theoremstyle{definition}

\newtheorem{defi}[theo]{Definition} \newtheorem{exam}[theo]{Example} 

\theoremstyle{remark}

\newtheorem{rema}[theo]{Remark}  

\newcommand {\absolute}[1] {\left| {#1} \right|}
\newcommand {\norm}[1] {\left\| {#1} \right\|}
\newcommand {\IGNORE}[1] {}

\begin{document}
\allowdisplaybreaks\frenchspacing

\setlength{\baselineskip}{1.1\baselineskip}

\title[Symbolic representations]{Symbolic representations of nonexpansive group automorphisms}

\author{Elon Lindenstrauss}

\address{Elon Lindenstrauss: Department of Mathematics, Princeton University, Princeton, NJ 08540, USA} \email{elonl@math.princeton.edu}

\author{Klaus Schmidt}

\address{Klaus Schmidt: Mathematics Institute, University of Vienna, Nordberg\-stra{\ss}e 15, A-1090 Vienna, Austria \newline\indent \textup{and} \newline\indent Erwin Schr\"odinger Institute for Mathematical Physics, Boltzmanngasse~9, A-1090 Vienna, Austria} \email{klaus.schmidt@univie.ac.at}

\begin{abstract}
If $\alpha $ is an irreducible nonexpansive ergodic automorphism of a compact abelian group $X$ (such as an irreducible nonhyperbolic ergodic toral automorphism), then $\alpha $ has no finite or infinite state Markov partitions, and there are no nontrivial continuous embeddings of Markov shifts in $X$. In spite of this we are able to construct a symbolic space $V$ and a class of shift-invariant probability measures on $V$ each of which corresponds to an $\alpha$-invariant probability measure on $X$. Moreover, every $\alpha$-invariant probability measure on $X$ arises essentially in this way.

The last part of the paper deals with the connection between the two-sided beta-shift $V_\beta $ arising from a Salem number $\beta $ and the nonhyperbolic ergodic toral automorphism $\alpha $ arising from the companion matrix of the minimal polynomial of $\beta $, and establishes an entropy-preserving correspondence between a class of shift-invariant probability measures on $V_\beta $ and certain $\alpha $-invariant probability measures on $X$. This correspondence is much weaker than, but still quite closely modelled on, the connection between the two-sided beta-shifts defined by Pisot numbers and the corresponding hyperbolic ergodic toral automorphisms.
\end{abstract}



\maketitle

\section{Introduction}

For expansive automorphisms $\alpha $ of compact connected abelian groups $X$, the attempt to find symbolic representations of the dynamical system $(X, \alpha )$ has a long and extensive history. In addition to the classical `geometric' constructions of Markov partitions (e.g. in \cite{AW}, \cite{Bowen} and \cite{Sinai}), there are explicit algebraic constructions of continuous equivariant finite-to-one maps from a sofic shift onto $X$. The first such construction for arbitrary irreducible hyperbolic toral automorphisms was given by R. Kenyon and A. Vershik in \cite{KV} (irreducibility is explained in Definition \ref{d:homoclinic}); a different, but related, general construction for irreducible expansive automorphisms of tori and solenoids was given by the second named author in \cite{S3}. In certain cases, this map can be chosen to be one-to-one almost everywhere (cf. \cite{S3} and \cite{SV}; of course, since $X$ is connected and a sofic shift completely disconnected one cannot hope to find a map which is one-to-one everywhere). The existence of such a map gives an explicit essentially one-to-one map between shift-invariant measures on a sofic shift and $\alpha $-invariant measures on $X$.

The key idea in the construction of these maps is to find a surjective equivariant map from some symbolic system $V$ onto $X$ (it turns out to be natural to set $V$ either equal to the space $\ell ^ \infty (\mathbb{Z} ,\mathbb{Z})$ of bounded integer sequences or to some sufficiently large compact shift-invariant subset of $\ell ^ \infty (\mathbb{Z} ,\mathbb{Z})$). Following an idea originally introduced by A. Vershik in \cite{Ver1}--\cite{Ver3} one may, for example, take a point $x \in X$ which is homoclinic to $0$ (i.e. which satisfies that $\lim_{|n|\to \infty }\alpha ^ nx=0$) and send any integer sequence $v=( \dots , v_{-1},v_0,v_1,\dots ) \in \ell ^ \infty (\mathbb{Z},\mathbb{Z})$ to the point
$$
\xi (v)=\sum_{ n \in \mathbb{Z} } v_n \alpha ^{ -n }x \in X.
$$
The resulting map $\xi \colon \ell ^ \infty (\mathbb{Z},\mathbb{Z})\longrightarrow X$ is equivariant (i.e. $\xi \circ \bar{\sigma }=\alpha \circ \xi $), and it is not hard to see that it is surjective. From this map $\xi $ one obtains a surjective map from the collection of shift-invariant probability measures on $\ell ^ \infty ( \mathbb{Z} , \mathbb{Z} )$ onto the $\alpha $-invariant measures on $X$. The more refined construction of \cite{S3} alluded to earlier is obtained by restricting this map to a carefully chosen sofic subshift $V \subset \ell ^ \infty ( \mathbb{Z} , \mathbb{Z} )$ on which $\xi $ is surjective and almost one-to-one. Other interesting and, indeed, more canonical examples arise when $\alpha $ is the automorphism of $X=\mathbb{T}^ m$ defined by the companion matrix of the minimal polynomial of a Pisot unit $\beta $ (i.e of an algebraic integer $\beta >1$ whose conjugates all have absolute values $<1$). In this case the corresponding two-sided beta-shift $V_\beta \subset \ell ^ \infty (\mathbb{Z},\mathbb{Z})$ is sofic, and the map $\xi \colon V_\beta \longrightarrow \mathbb{T}^ m$ defined above is surjective, finite-to-one and conjectured to be almost one-to-one (cf. \cite{S3}--\cite{SV}).

In this paper we investigate to what extent one can find a suitable substitute for this construction in the nonexpansive case. This question is motivated by the somewhat exotic behaviour of invariant probability measures of irreducible nonhyperbolic ergodic toral automorphisms described in \cite{LS1}: if $\mu $ is a probability measure on $X=\mathbb{T}^ n$ which is invariant under an irreducible nonhyperbolic ergodic toral automorphism $\alpha $, but which is completely singular with respect to Lebesgue measure, then there exists an $\alpha $-invariant Borel set $B \subset X$ which intersects $\mu \textsl{-a.e.}$ coset of the dense central subgroup $X ^{(0)}\subset X$, on which $\alpha $ acts isometrically, in at most one point. If $\mu $ is weakly mixing then one may assume in addition that $\mu (B)=1$. Any natural `symbolic model' of such an automorphism would enable one to construct such measures quite explicitly.

The first difficulty one encounters in the search for symbolic models of an irreducible ergodic nonexpansive automorphism $\alpha $ of a compact connected abelian group $X$ is that every continuous equivariant map $\phi \colon Y \longrightarrow X$ from a mixing shift of finite type $Y$ with finite or countably infinite alphabet (or from a two-sided beta-shift) to $X$ maps the shift space to a single point (cf. Corollaries \ref{c:xi*}--\ref{c:xi*4} and Remark \ref{r:xi*}); in particular, $(X,\alpha )$ cannot have finite or countably infinite Markov partitions in any reasonable sense. The reason for this is that these automorphisms have no nonzero homoclinic points (Theorem \ref{t:1}).

It is, however, possible to define a continuous map $\tilde \xi$ from the noncompact space $\tilde Y = \ell ^ \infty (\mathbb{Z},\mathbb{Z}) \times X ^{(0)}$ to $X$ which is equivariant with respect to an isometric cocycle extension $\tilde{\sigma }\colon \tilde{Y}\longrightarrow \tilde{Y}$ of the shift on $\ell ^ \infty (\mathbb{Z},\mathbb{Z})$. This map is surjective (though far from injective), and allows us in particular to map shift-invariant probability measures on $\tilde Y$ to $\alpha$-invariant measures on $X$. Indeed, we show the following (cf. Proposition \ref{p:measures}; \emph{central equivalence} is explained in Definition \ref{d:equivalent}).

\begin{theo}
\label{many to one theorem}
Any $\alpha$-invariant probability measure on $X$ is centrally equivalent to the push-forward under $\tilde \xi$ of a $\tilde{\sigma }$-invariant probability measure on $ \tilde Y$, which may further be taken to be compactly supported.
\end{theo}

We emphasize that this is true even for Lebesgue (or Haar) measure (since central equivalence preserves entropy, any $\alpha $-invariant probability measure on $X$ which is centrally equivalent to Lebesgue measure must be equal to Lebesgue measure). Since $\tilde Y$ is a noncompact extension of $\ell ^ \infty (\mathbb{Z},\mathbb{Z})$, not every shift invariant measure on $\ell ^ \infty (\mathbb{Z},\mathbb{Z})$ can be lifted to a $\tilde{\sigma }$-invariant probability measure on $ \tilde Y$. The measures which \emph{can} be lifted are precisely those for which the cocycle appearing in the definition of $\tilde{\sigma }$ is a coboundary (cf. Theorem \ref{t:3} and Proposition \ref{p:measures}). The fact that it is natural to consider only those measures on a symbolic model for which this cocycle is trivial can be viewed as a manifestation of some weak form of measure rigidity for nonexpansive group automorphisms.

The main drawback of Theorem~\ref{many to one theorem} is that the same measure on $X$ can be obtained as the push-forward of many measures on $\tilde Y$; furthermore, it is quite hard to understand properties such as the entropy of the resulting measures in terms of the properties of the original measure. In order to resolve such difficulties one would like to replace $\ell ^ \infty (\mathbb{Z},\mathbb{Z})$ by a smaller closed subshift, just like in the hyperbolic case.

In the case of toral automorphisms corresponding to Pisot numbers (i.e. of irreducible hyperbolic toral automorphisms with one-dimensional unstable manifolds) there is a natural candidate: the beta-shift $V_\beta$ corresponding to the unique `large' eigenvalue $\beta$ of the automorphism. Motivated by this question we devote Section \ref{s:symbolic} to a problem which has also provided much of the original motivation for this research: the connection between the two-sided beta-shift $V_\beta $ arising from a Salem number $\beta $ and the nonhyperbolic ergodic toral automorphism $\alpha $ defined by the companion matrix of the minimal polynomial of $\beta $ (a \emph{Salem number} is an algebraic unit $\beta >1$ whose conjugates all have absolute values $\le1$, with at least one conjugate of absolute value $=1$). In contrast to the Pisot case, which is reasonably well understood (though some important questions in this construction are still unresolved, as described \vpageref{one-to-one}), the beta-shifts associated with Salem numbers still hold many mysteries. For example, it is not known whether they are always sofic (cf. \cite{Boyd1}--\cite{Boyd3} and \cite{Sbeta}). Not surprisingly, the dynamical interpretation of two-sided beta-shifts arising from Salem numbers is much more complicated than in the Pisot case.

By restricting the map $\tilde \xi\colon \tilde{Y}\longrightarrow X$ described above to the space $\tilde Y _ \beta = V _ \beta \times X ^ {( 0 )}$ we obtain a map from $\tilde{\sigma }$-invariant probability measures on $\tilde Y _ \beta$ (or, equivalently, from shift invariant probability measures on $V _ \beta$ satisfying the cocycle condition mentioned above) to $\alpha$-invariant measures on $X$. In particular, the following theorem follows from the main result of Section \ref{s:symbolic} (cf. Theorem \ref{t:weaklybounded}).

\begin{theo}
\label{theorem: same entropy}
For any $\tilde{\sigma }$-invariant probability measure $\tilde \mu$ on $\tilde Y _ \beta$, the entropy of the push-forward $ \tilde \xi _ {*} (\tilde \mu)$ is equal to that of $\tilde \mu$.
% at If "squiggle Greek mu" is compactly supported, then " support of squiggle Greek xi _ {*} (squiggle Greek mu)" is a proper subset of "X" and "squiggle Greek xi|_{support of squiggle Greek mu}" is bounded to one.
\end{theo}

By constructing in Section \ref{s: examples} shift-invariant probability measures on $V _ \beta$ satisfying a strong form of the cocycle condition with entropies arbitrarily close to $\log \beta$ we obtain from Theorem~\ref{theorem: same entropy} $\alpha$-invariant probability measures on $X$ which are singular with respect to Lebesgue measure and whose entropies are arbitrarily close to $\log\,\beta $, the topological entropy of $(X, \alpha)$.

%%%%CHANGES
In the course of proving  of Theorem~\ref{theorem: same entropy} , we show that
Lebesgue measure on $X$ cannot be represented as $\tilde \xi _ {*} (\tilde \mu)$ with $\tilde \mu$ a measure on $\tilde Y _ \beta$ as above. The main question highlighted by our work is the following:

\begin{ques}
Can every $\alpha$-invariant probability measure on $X$ which is completely singular with respect to Lebesgue measure be presented as $\tilde \xi _ {*} (\tilde \mu)$ for an invariant probability measure $\tilde \mu$ on $\tilde Y _ \beta$?
\end{ques}

At present, we have no evidence in either direction. Even if the answer turns out to be negative, it would be interesting to understand the relation between the space of measures obtained by the construction of Theorem~\ref{theorem: same entropy} and the space of all invariant measures.

It follows from Theorem~\ref{theorem: same entropy} and the results of Section \ref{s: examples} that $\tilde \xi _ {*} (\tilde Y _ \beta) \subset X$ is fairly large; for example, it can be shown that its Hausdorff dimension is the same as that of $X$. However, we do not even know the answer to the following natural question:

\begin{ques}
\label{question about pseudo-covers}
Is $\tilde \xi _ {*} (\tilde Y _ \beta) = X$?
\end{ques}

In the notation of Section~\ref{s:nonexpansive}, Question \ref{question about pseudo-covers} can be rephrased as follows: is the $\beta$-shift $V _ \beta$ a pseudo-cover of $X$?

We end this introduction with a comment on a technical simplification we adopt throughout this paper: every irreducible automorphism $\alpha $ of a compact connected abelian group $X$ is finitely equivalent to a group automorphism of the special form $\alpha _{R_1/(f)}$ described in \eqref{eq:principal}--\eqref{eq:alpha2}. By restricting ourselves to automorphisms of this special form we avoid some minor notational and technical complications in the statements of our results due to the presence of finite-to-one factor maps, but our discussion here can be translated to the general case without any difficulty.

\section{Homoclinic points of irreducible group automorphisms}
\label{s:homoclinic}

\begin{defi}
\label{d:homoclinic}
Let $\alpha $ be a continuous automorphism of a compact abelian group $X$ with identity element $0=0_X$. The automorphism $\alpha $ is \emph{irreducible} if every closed $\alpha $-invariant subgroup $Y \subsetneq X$ is finite. A point $x \in X$ is \emph{$\alpha $-homoclinic} (or simply \emph{homoclinic}) if $\lim_{|n|\to \infty }\alpha ^ nx=0$. The set $\Delta _\alpha (X)$ of homoclinic points in $X$ is an $\alpha $-invariant subgroup.
\end{defi}

Recall that two continuous automorphisms $\alpha $ and $\beta $ of compact abelian groups $X$ and $Y$ are \emph{finitely equivalent} if there exist continuous, surjective, equivariant and finite-to-one group homomorphisms $\phi \colon X \longrightarrow Y$ and $\psi \colon Y \longrightarrow X$. In order to describe all irreducible automorphisms of compact abelian groups up to finite equivalence we use notation from \cite{LS1}. Let $R_1=\mathbb{Z}[u ^{\pm 1}]$ be the ring of Laurent polynomials with integral coefficients. Every $h \in R_1$ is of the form
\begin{equation}
\label{eq:h}
h=\sum_{k \in \mathbb{Z}}h_ku ^ k
\end{equation}
with $h_k \in \mathbb{Z}$ for every $k \in \mathbb{Z}$ and $h_k=0$ for all but finitely many $k$. Fix an irreducible polynomial
\begin{equation}
\label{eq:f}
f=f_0+\dots +f_mu ^ m \in R_1
\end{equation}
with $m>0$, $f_m>0$ and $f_0\ne0$, denote by $\Omega _f$ the set of roots of $f$, and set
\begin{equation}
\label{eq:Omega}
\begin{gathered}
\Omega _f ^-=\{ \omega \in \Omega _f:|\omega |<1 \},\enspace \Omega _f ^{(0)} =\{ \omega \in \Omega _f:|\omega |=1 \},
\\
\Omega _f ^+=\{ \omega \in \Omega _f:|\omega |>1 \}.
\end{gathered}
\end{equation}
We write $\mathbb{T}=\mathbb{R}/\mathbb{Z}$ for the circle group, define the shift $\tau \colon \mathbb{T}^ \mathbb{Z}\longrightarrow \mathbb{T}^ \mathbb{Z}$ by
\begin{equation}
\label{eq:tau}
\tau (x)_n=x_{n+1}
\end{equation}
for every $x=(x_n)\in \mathbb{T}^ \mathbb{Z}$, and set
\begin{equation}
\label{eq:h of tau}
h(\tau )=\sum_{k \in \mathbb{Z}}h_k \tau ^ k\colon \mathbb{T}^ \mathbb{Z}\longrightarrow \mathbb{T}^ \mathbb{Z}
\end{equation}
for every $h \in R_1$ of the form \eqref{eq:h}. Consider the closed, shift-invariant subgroup
\begin{align}
X=X_{R_1/(f)}&=\biggl\{ x \in \mathbb{T}^{\mathbb{Z}}: \sum\nolimits_{n \in \mathbb{Z}} f_nx_{k+n}=0 \pmod 1\enspace \text{for every}\enspace k \in \mathbb{Z}\biggr\}\notag
\\
&=\{ x \in \mathbb{T}^{\mathbb{Z}}: f(\tau )(x)=0 \} = \ker f(\tau ),
\label{eq:principal}
\end{align}
and write
\begin{equation}
\label{eq:alpha2}
\alpha =\alpha _{R_1/(f)}
\end{equation}
for the restriction of $\tau $ to $X \subset \mathbb{T}^{\mathbb{Z}}$ (cf. \cite[(2.3) and (2.11)]{LS1}). By \cite[Theorem 7.1 and Propositions 7.2--7.3]{DSAO}, $\alpha $ is nonexpansive if and only if $\Omega _f ^{(0)} \ne\varnothing $, and ergodic if and only if $f$ is not cyclotomic (i.e. if and only if $f$ does not divide $u ^ m-1$ for any $m\ge1$). In view of this we adopt the following terminology.
\begin{defi}
\label{d:hyperbolic}
The polynomial $f$ in \eqref{eq:f} is \emph{hyperbolic} if $\Omega _f ^{(0)} =\varnothing $, \emph{nonhyperbolic} if $\Omega _f ^{(0)} \ne\varnothing $, and \emph{cyclotomic} if $\Omega _f ^{(0)}$ contains a root of unity.
\end{defi}

According to \cite{S2}, every irreducible automorphism $\alpha $ of a compact abelian group $X$ is finitely equivalent to an automorphism of the form $\alpha _{R_1/(f)}$ for some irreducible polynomial $f \in R_1$. Note that the automorphisms $\alpha $ and $\alpha _{R_1/(f)}$ are expansive if and only if $f$ is hyperbolic, and ergodic if and only if $f$ is not cyclotomic.

\emph{For the remainder of this article we assume that the irreducible polynomial $f$ in \eqref{eq:f} is noncyclotomic.}

We denote by $\| \cdot \|_1$ and $\| \cdot \|_\infty $ the norms on the Banach spaces $\ell ^ 1(\mathbb{Z},\mathbb{R})$ and $\ell ^ \infty (\mathbb{Z},\mathbb{R})$ and write $\ell ^ 1(\mathbb{Z},\mathbb{Z}) \subset \ell ^ 1(\mathbb{Z},\mathbb{R})$ and $\ell ^ \infty (\mathbb{Z},\mathbb{Z}) \subset \ell ^ \infty (\mathbb{Z},\mathbb{R})$ for the subgroups of integer-valued functions. By viewing every $h=\sum_{n \in \mathbb{Z}}h_nu ^ n \in R_1$ as the element $(h_n)\in \ell ^ 1(\mathbb{Z},\mathbb{Z})$ we can identify $R_1$ with $\ell ^ 1(\mathbb{Z},\mathbb{Z})$.

We furnish the space $\ell ^ \infty (\mathbb{Z},\mathbb{R})$ with the topology of coordinate-wise convergence. In this topology $\ell ^ \infty (\mathbb{Z},\mathbb{R})$ is a metrizable topological group, $\ell ^ \infty (\mathbb{Z},\linebreak[0]\mathbb{Z})\subset \ell ^ \infty (\mathbb{Z},\mathbb{R})$ is a closed subgroup, and the shift-invariant sets
\begin{equation}
\label{eq:Br}
\begin{gathered}
B_r(\ell ^ \infty (\mathbb{Z},\mathbb{R}))=\{ w \in \ell ^ \infty (\mathbb{Z},\mathbb{R}):\| w \|_\infty \le r \},
\\
B_r(\ell ^ \infty (\mathbb{Z},\mathbb{Z}))=B_r(\ell ^ \infty (\mathbb{Z},\mathbb{R}))\cap \ell ^ \infty (\mathbb{Z},\mathbb{Z})
\end{gathered}
\end{equation}
are compact for every $r\ge0$: on these sets our topology coincides with the weak${}^*$-topology.

As in \cite{ES} we denote by $\bar{\sigma }$ the shift
\begin{equation}
\label{eq:sigmabar}
(\bar \sigma w)_n=w_{n+1}
\end{equation}
on $\ell ^ \infty (\mathbb{Z},\mathbb{R})$, observe that $\bar{\sigma }\colon \ell ^ \infty (\mathbb{Z},\mathbb{R})\longrightarrow \ell ^ \infty (\mathbb{Z},\mathbb{R})$ is a continuous group automorphism, and define, for every $h=\sum_{k \in \mathbb{Z}}h_ku ^ k \in R_1$, a continuous group homomorphism
\begin{equation}
\label{eq:hsigmabar}
h(\bar \sigma )=\sum_{k \in \mathbb{Z}}h_k \bar \sigma ^ k\colon \ell ^ \infty (\mathbb{Z},\mathbb{R})\longrightarrow \ell ^ \infty (\mathbb{Z},\mathbb{R})
\end{equation}
(cf. \eqref{eq:h of tau}). The map $\rho \colon \ell ^ \infty (\mathbb{Z},\mathbb{R})\longrightarrow \mathbb{T}^ \mathbb{Z}$, given by
\begin{equation}
\label{eq:rho}
\rho (w)_n=w_n\;(\textup{mod}\;1)
\end{equation}
for every $w=(w_n)\in \ell ^ \infty (\mathbb{Z},\mathbb{R})$ and $n \in \mathbb{Z}$, is a continuous surjective group homomorphism with
\begin{equation}
\label{eq:equivariant}
\rho \circ \bar{\sigma }=\tau \circ \rho ,
\end{equation}
and the set
\begin{equation}
\label{eq:Wf}
W_f=\rho ^{-1}(X)=f(\bar{\sigma })^{-1}(\ell ^ \infty (\mathbb{Z},\mathbb{Z}))\subset \ell ^ \infty (\mathbb{Z},\mathbb{R}),
\end{equation}
is a closed and shift-invariant subgroup with $\ker \rho =\ell ^ \infty (\mathbb{Z},\mathbb{Z})\subset W_f$.

The kernel
\begin{equation}
\label{eq:Wf0}
W_f ^{(0)} =\ker f(\bar{\sigma })\subset W_f
\end{equation}
is obviously finite-dimensional and the restriction of $\bar{\sigma }$ to the complexification $\mathbb{C} \otimes_\mathbb{R}W_f ^{(0)} $ of $W_f ^{(0)} $ is linear. Hence $\bar{\sigma }$ has a nonzero eigenvector $v \in \mathbb{C}\otimes_\mathbb{R}W_f ^{(0)} $ with eigenvalue $\omega \in \mathbb{C}$, say, and $f(\omega )=0$. As $\bar{\sigma }$ is an isometry on $W_f ^{(0)} $ we conclude that $\omega \in \Omega _f ^{(0)} $. Conversely, if $\omega \in \Omega _f ^{(0)} $, then we set $v_n=\omega ^ n$ for every $n \in \mathbb{Z}$ and obtain that $v=(v_n)\in \mathbb{C} \otimes_\mathbb{R}W_f ^{(0)} $.

This shows that $W_f ^{(0)} =\ker f(\bar{\sigma })\subset W_f \subset \ell ^ \infty (\mathbb{Z},\mathbb{R})$ is the linear span of the vectors $\{ \Re (w(\omega )),\Im (w(\omega )):\omega \in \Omega _f ^{(0)} \}$ with
\begin{equation}
\label{eq:womega}
w(\omega )_n=\omega ^ n, \qquad \Re (w(\omega ))_n=\Re(\omega ^ n), \qquad \Im (w(\omega ))_n = \Im (\omega ^ n)
\end{equation}
for every $n \in \mathbb{Z}$ and $\omega \in \Omega _f ^{(0)} $, where $\Re$ and $\Im$ denote the real and imaginary parts. By \eqref{eq:equivariant},
\begin{equation}
\label{eq:X0}
X ^{(0)}=\rho (\ker f(\bar{\sigma }))=\rho (W_f ^{(0)} )
\end{equation}
is an $\alpha $-invariant subgroup of $X$, and the irreducibility of $\alpha $ implies that the closure of $X ^{(0)}$ is either equal to $\{ 0 \}$ (if $\alpha $ is expansive), or to $X$ (if $\alpha $ is nonexpansive). The group $X ^{(0)}\subset X$ in \eqref{eq:X0} is isomorphic to $W_f ^{(0)}$, since $\rho $ is injective on $W_f ^{(0)}$, and coincides with the central subgroup of $X$ defined in \cite[(3.3)]{LS1} on which $\alpha $ acts isometrically.

We write
$$
\smash[t]{\frac 1{f(u)}= \frac 1{f_m}\sum_{\omega \in \Omega _f}\frac{b_\omega }{u-\omega }}
$$
for the partial fraction decomposition of $1/f$ with $b_\omega \in \mathbb{C}$ for every $\omega \in \Omega _f$ and define elements $w ^{\Delta _\pm}$ and $w ^{\Delta _0}$ in $\ell ^ \infty (\mathbb{Z},\mathbb{R})$ by
\begin{equation}
\label{eq:homoclinic}
\begin{aligned}
w ^{\Delta _+}_n&=
\begin{cases}
\frac 1{f_m}\cdot \sum_{\omega \in \Omega _f ^-}\hspace{10mm}b_\omega \omega ^{n-1}&\hspace{5mm}\textup{if}\enspace n\ge1,
\\
\frac 1{f_m}\cdot \sum_{\omega \in \Omega _f ^{(0)} \cup \Omega _f ^+}\hspace{1.2mm}-b_\omega \omega ^{n-1}&\hspace{5mm}\textup{if}\enspace n\le0,
\end{cases}
\\
w ^{\Delta _-}_n&=
\begin{cases}
\frac 1{f_m}\cdot \sum_{\omega \in \Omega _f ^-\cup \Omega _f ^{(0)} }\hspace{4.1mm}b_\omega \omega _i ^{n-1}&\hspace{5mm}\textup{if}\enspace n\ge1,
\\
\frac 1{f_m}\cdot \sum_{\omega \in \Omega _f ^+}\hspace{7mm}-b_\omega \omega ^{n-1}&\hspace{5mm}\textup{if}\enspace n\le0,
\end{cases}
\\
w ^{\Delta _0}_n&=\hspace{3.1mm}\tfrac 1{f_m}\cdot \textstyle \sum_{\omega \in \Omega _f ^{(0)} }\hspace{10.3mm}b_\omega \omega ^{n-1}\hspace{9mm}\textup{for every}\enspace n \in \mathbb{Z}.
\end{aligned}
\end{equation}
Then
\begin{equation}
\label{eq:homoclinic2}
\begin{gathered}
w ^{\Delta _0}\in W_f ^{(0)} ,\qquad w ^{\Delta _+}+w ^{\Delta _0}=w ^{\Delta _-},
\\
f(\bar{\sigma })(w ^{\Delta _+})_n=f(\bar{\sigma })(w ^{\Delta _-})_n=v ^ \Delta _n=
\begin{cases}
1&\textup{if}\enspace n=0,
\\
0&\textup{otherwise},
\end{cases}
\end{gathered}
\end{equation}
where we are using the formal power series identities
\begin{align*}
\sum_{n \in \mathbb{Z}}w ^{\Delta _+}_nu ^ n&= \frac 1{f_m}\cdot \biggl(\sum_{\omega \in \Omega _f ^-}\frac{b_\omega u}{1-\omega u}+\sum_{\omega \in \Omega _f ^{(0)} \cup \Omega _f ^+}\frac{-b_\omega \omega ^{-1}}{1-\omega ^{-1}u ^{-1}}\biggr)
\\
&=\frac 1{f_m}\cdot \sum_{\omega \in \Omega _f} \frac{b_\omega }{u ^{-1}-\omega }=\frac 1{f(u ^{-1})}
\\
&=\smash[b]{\frac 1{f_m}\cdot \biggl(\sum_{\omega \in \Omega _f ^-\cup \Omega _f ^{(0)} } \frac{b_\omega u}{1-\omega u}+\sum_{\omega \in \Omega _f ^+}\frac{-b_\omega \omega ^{-1}}{1-\omega ^{-1}u ^{-1}}\biggr)=\sum_{n \in \mathbb{Z}}w ^{\Delta _-}_nu ^ n}
\end{align*}
and
$$
\smash[t]{\sum_{n \in \mathbb{Z}}f(\bar{\sigma })(w)_nu ^ n=f(u ^{-1})\cdot \sum_{n \in \mathbb{Z}}w_nu ^ n}
$$
for every $w=(w_n)\in \ell ^ \infty (\mathbb{Z},\mathbb{R})$. The points $w ^{\Delta _\pm}\in \ell ^ \infty (\mathbb{Z},\mathbb{R})$ have the following properties.
\begin{align}
&x ^{\Delta _\pm}=\rho (w ^{\Delta _\pm})\in X&\textup{by \eqref{eq:homoclinic2}},\notag
\\
&\hspace{-1mm}\lim_{n \to \infty }w ^{\Delta _+}_n=\lim_{n \to \infty }w ^{\Delta _-}_{-n}=0&\textup{exponentially fast},
\label{eq:decay}
\\
&x ^{\Delta _+}=x ^{\Delta _-}&\textup{if and only if $\alpha $ is expansive}.\notag
\end{align}

\section{A review of the expansive case}
\label{s:expansive}

One of the key tools in attempting to find symbolic covers or representations of the automorphism $\alpha =\alpha _{R_1/(f)}$ lies in identifying the subgroup
\begin{equation}
\label{eq:Vf}
V_f=f(\bar{\sigma })(W_f)\subset \ell ^ \infty (\mathbb{Z},\mathbb{Z}).
\end{equation}
We first discuss the space $V_f$ in the expansive setting, before moving on to the nonexpansive case.

Suppose that the polynomial $f$ in \eqref{eq:f} is hyperbolic (i.e. that $\Omega _f ^{(0)} =\varnothing $ in \eqref{eq:Omega}). In this case
\begin{equation}
\label{eq:fundamental}
w ^{\Delta _+}=w ^{\Delta _-}=w ^ \Delta ,\qquad x ^{\Delta _+}=x ^{\Delta _-}=x ^ \Delta ,\qquad w ^{\Delta _0}=0.
\end{equation}
The point $x ^ \Delta $ is a \emph{fundamental homoclinic point} of $\alpha $ in the sense of \cite{LS}:
\begin{equation}
\label{eq:fundamental2}
\Delta _\alpha (X)=\{ h(\alpha )(x ^ \Delta ):h \in R_1 \}.
\end{equation}

In the case where $f_m=|f_0|=1$ in \eqref{eq:f} and $X$ is therefore isomorphic to $\mathbb{T}^ m=\mathbb{R}^ m/\mathbb{Z}^ m$, the fundamental homoclinic point $x ^ \Delta $ has a convenient geometric description.
\label{description1}
The automorphism $\alpha =\alpha _{R_1/(f)}$ of $X=X_{R_1/(f)}\subset \mathbb{T}^ \mathbb{Z}$ in \eqref{eq:principal} is algebraically conjugate to the companion matrix
\begin{equation}
\label{eq:companion}
M_f= \left[
\begin{smallmatrix}
0&1&0&\dots&0&0
\\
0&0&1&\dots&0&0
\\
\vdots&&\vdots&\ddots&\vdots&0
\\
0&0&0&\dots&0&1
\\
-f_0&-f_1&-f_2&\dots&-f_{m-2}&-f_{m-1}
\end{smallmatrix}
\right]\!\!,
\end{equation}
of $f$, acting on $\mathbb{T}^{m}$ from the left, where the isomorphism between $X$ and $\mathbb{T}^{m}$ is the coordinate projection
$$
\smash[t]{x \mapsto \left[
\begin{smallmatrix}
x_0
\\
x_1
\\
\vdots
\\
x_{m-1}
\end{smallmatrix}
\right]\!\!.}
$$
We write $W ^{(s)}\subset \mathbb{R}^ m$ and $W ^{(u)}\subset \mathbb{R}^ m$ for the contracting and expanding subspaces of the matrix $M_f$. The quotient map $\pi \colon \mathbb{R}^ m \longrightarrow \mathbb{T}^ m$ is injective on $W ^{(s)}$ and $W ^{(u)}$, and the dense subgroups $X ^{(s)}=\pi (W ^{(s)})$ and $X ^{(u)}=\pi (W ^{(u)})$ satisfy that
$$
\Delta _\alpha (X)= X ^{(s)}\cap X ^{(u)}=\pi ((W ^{(s)}+\mathbb{Z}^ m)\cap W ^{(u)}).
$$
There exists a unique point $y ^ \Delta \in (W ^{(s)}+\mathbf{e}^{(1)})\cap W ^{(u)}$, where $\mathbf{e}^{(1)}=(1,0,\dots )$ is the first unit vector in $\mathbb{R}^ m$. Since $M_f$ is of the form \eqref{eq:companion}, the orbit $\{ M_f ^ n \mathbf{e}^{(1)}:n \in \mathbb{Z}\} \subset \mathbb{Z}^ m$ generates $\mathbb{Z}^ m$ as a group, which is easily seen to imply that the homoclinic point $x ^ \Delta =\pi (y ^ \Delta )$ is indeed fundamental.

We return to our more general setting. From \eqref{eq:decay}--\eqref{eq:fundamental} it follows that
$$
\smash{\| w ^ \Delta \|_1=\sum_{n \in \mathbb{Z}}|w_n ^ \Delta |<\infty ,}
$$
and that
$$
\smash[t]{\bar \xi (v)=\sum_{n \in \mathbb{Z}}v_n \bar \sigma ^{-n}w ^ \Delta }
$$
is a well-defined element of $\ell ^ \infty (\mathbb{Z},\mathbb{R})$ for every $v \in \ell ^ \infty (\mathbb{Z},\mathbb{Z})$. As in \cite{ES} we denote by
\begin{equation}
\label{eq:211}
\bar \xi \colon \ell ^ \infty (\mathbb{Z},\mathbb{Z})\longrightarrow \ell ^ \infty (\mathbb{Z},\mathbb{R}),\qquad \xi =\rho \circ \bar \xi \colon \ell ^ \infty (\mathbb{Z},\mathbb{Z})\longrightarrow X
\end{equation}
the resulting continuous group homomorphisms and observe that

\begin{equation}
\label{eq:xi}
\xi (v)=\sum_{n \in \mathbb{Z}}v_n \alpha ^{-n}x ^ \Delta
\end{equation}
for every $v \in \ell ^ \infty (\mathbb{Z},\mathbb{Z})$. Hence
$$
\xi \circ \bar \sigma =\alpha \circ \xi ,
$$
i.e. $\xi $ is equivariant. We summarize this discussion in a theorem; the relevant proofs can be found in \cite{ES}.

\begin{theo}
\label{t:xi}
Let $f \in R_1$ be an irreducible hyperbolic polynomial, and let $\alpha =\alpha _{R_1/(f)}$ be the expansive automorphism of the compact connected abelian group $X=X_{R_1/(f)}$ defined in \eqref{eq:principal}--\eqref{eq:alpha2}. Then
$$
V_f=f(\bar{\sigma })(W_f)=\ell ^ \infty (\mathbb{Z},\mathbb{Z}),
$$
and the homomorphisms $\bar \xi \colon \ell ^ \infty (\mathbb{Z},\mathbb{Z})\longrightarrow \ell ^ \infty (\mathbb{Z},\mathbb{R})$ and $\xi =\rho \circ \bar \xi \colon \ell ^ \infty (\mathbb{Z},\mathbb{Z})\linebreak[0]\longrightarrow X$ in \eqref{eq:211}--\eqref{eq:xi} satisfy that
\begin{gather*}
\xi (\ell ^ \infty (\mathbb{Z},\mathbb{Z}))=X,
\\
\ker \xi =f(\bar \sigma )(\ell ^ \infty (\mathbb{Z},\mathbb{Z}))\subset \ell ^ \infty (\mathbb{Z},\mathbb{Z}),
\\
\xi \circ \bar{\sigma }=\alpha \circ \xi .
\end{gather*}
\end{theo}

In \cite{S3} it was proved that there always exists a compact shift-invariant subset (in fact, a sofic subshift) $\tilde{V}\subset \ell ^ \infty (\mathbb{Z},\mathbb{Z})$ such that the restriction of $\xi $ to $\tilde{V}$ is surjective and almost one-to-one (for the definition of a sofic shift we refer to \cite{LM} and \cite{Weiss}). In general, however, there is at present no \emph{distinguished} candidate for such a set $\tilde{V}$.

In Section~\ref{s: beta} we present an interesting special case in which there is a natural candidate for $\tilde{V}$ (cf. \cite{S3}--\cite{Ver3}): the beta-shift.

\section{Homoclinic points and coding in the nonexpansive case}
\label{s:nonexpansive}

Now suppose that the irreducible polynomial $f$ in \eqref{eq:f} is nonhyperbolic and not cyclotomic, and that the ergodic automorphism $\alpha =\alpha _{R_1/(f)}$ of the compact connected abelian group $X=X_{R_1/(f)}$ is therefore ergodic and nonexpansive. Since $f$ is irreducible and has a root of absolute value $1$, $m$ is even and $f_i=f_{m-i}$ for $i=0,\dots ,m$, and we assume that $f_0=f_m>0$. In contrast to the expansive situation, $W_f ^{(0)} =\ker f(\bar{\sigma })\ne \{ 0 \}$, the central subgroup $X ^{(0)}=\rho (W_f ^{(0)} )$ in \eqref{eq:X0} is dense in $X$, and the following theorem shows that there are no nonzero $\alpha $-homoclinic points.

\begin{theo}
\label{t:1}
Let $f \in R_1$ be an irreducible nonhyperbolic polynomial which is not cyclotomic, and let $\alpha =\alpha _{R_1/(f)}$ be the ergodic and nonexpansive automorphism of the compact connected abelian group $X=X_{R_1/(f)}$ defined in \eqref{eq:principal}--\eqref{eq:alpha2}. Then $\Delta _\alpha (X)=\{ 0 \}$.
\end{theo}

\begin{coro}
\label{c:1}
Let $\alpha $ be an irreducible, ergodic and nonexpansive automorphism of a compact connected abelian group $X$. Then $\Delta _\alpha (X)=\{ 0 \}$.
\end{coro}

\begin{proof}[Proof of Theorem \ref{t:1}]
The triviality of the homoclinic group $\Delta _\alpha (X)$ for irreducible nonhyperbolic ergodic toral automorphisms was shown in \cite{LS}. Here we give another (and slightly more general) proof using the methods described in the previous section.

Suppose that $x$ is a nonzero $\alpha $-homoclinic point. Since the restriction of $\alpha $ to $X ^{(0)}$ is an isometry it is clear that $X ^{(0)}\cap \Delta _\alpha (X)=\{ 0 \}$ and hence that $x\notin X ^{(0)}$.

We choose $w \in W_f \subset \ell ^ \infty (\mathbb{Z},\mathbb{R})$ such that $\rho (w)=x$ and $\lim_{|n|\to \infty }w_n=0$ (such a choice is obviously possible). Then $v=f(\bar{\sigma })(w)\in \ell ^ \infty (\mathbb{Z},\mathbb{Z})$ has only finitely many nonzero coordinates and is therefore of the form $h(\bar{\sigma })(v ^ \Delta )$ for some $h \in R_1$, where the point $v ^ \Delta \in \ell ^ \infty (\mathbb{Z},\mathbb{Z})$ is defined in \eqref{eq:homoclinic2}. We put
$$
w ^*=h(\bar{\sigma })(w ^{\Delta _-})\in \ell ^ \infty (\mathbb{Z},\mathbb{R})
$$
and observe that
$$
w-w ^*\in \ker(f(\bar{\sigma }))=W_f ^{(0)}
$$
by \eqref{eq:homoclinic2}. From \eqref{eq:homoclinic} we know that $\lim_{n \to-\infty }w ^*_n=\lim_{n \to-\infty }w_n=0$, and hence that $w=w ^*$ and $\lim_{n \to \infty }w_n ^*=0$, since every element in $W_f ^{(0)}$ is almost periodic. However,
$$
w_n ^*=\frac 1{f_m}\sum_{\omega \in \Omega _f ^{(0)} \cup \Omega _f ^+}b_\omega \omega ^{n-1}h(\omega )
$$
for all sufficiently large positive $n$, which shows that
\begin{equation}
\label{eq:sum}
\sum_{\omega \in \Omega _f ^{(0)} }b_\omega \omega ^ nh(\omega )=0\enspace \textup{for every}\enspace n\ge0.
\end{equation}
From \eqref{eq:sum} we see that $h(\omega )=0$ for every $\omega \in \Omega _f ^{(0)} $ or, equivalently, that $h$ is divisible by $f$. We set $h=fh'$ with $h'\in R_1$, $v'=h'(\bar{\sigma })(v ^ \Delta )\in \ell ^ \infty (\mathbb{Z},\mathbb{Z})$ and $w'=h'(\bar{\sigma })(w ^{\Delta _-})\in W_f$ as above, and conclude that $w=w ^*=f(\bar{\sigma })(w')\in \ell ^ \infty (\mathbb{Z},\mathbb{Z})$ and $x=\rho (w)=0$, contrary to our choice of $x$.
\end{proof}

\begin{proof}[Proof of Corollary \ref{c:1}]
There exists an irreducible, nonhyperbolic and noncyclotomic polynomial $f \in R_1$ such that $\alpha $ is finitely equivalent to $\alpha _{R_1/(f)}$ (cf. \eqref{eq:principal}--\eqref{eq:alpha2}). If $\phi \colon X \longrightarrow X_{R_1/(f)}$ is a continuous, finite-to-one and equivariant group homomorphism, then the restriction of $\phi $ to $\Delta _\alpha (X)$ is injective and $\phi (\Delta _\alpha (X))\subset \Delta _{\alpha _{R_1/(f)}}(X_{R_1/(f)})=\{ 0 \}$ by Theorem \ref{t:1}. This proves that $\Delta _\alpha (X)=\{ 0 \}$.
\end{proof}

Although $\alpha =\alpha _{R_1/(f)}$ has no nonzero homoclinic points, we have at our disposal the `one-sided homoclinic' points $x ^{\Delta _\pm}$ in \eqref{eq:decay}. Again it may be helpful to identify these points in the special case where $f_0=f_m=1$ in \eqref{eq:f} and $X$ is therefore isomorphic to $\mathbb{T}^ m=\mathbb{R}^ m/\mathbb{Z}^ m$ (remember that $f(u)=u ^ mf(u ^{-1})$!). As \vpageref{description1} we write $W ^{(s)}\subset \mathbb{R}^ m, W ^{(u)}\subset \mathbb{R}^ m$ and $W ^{(0)}\subset \mathbb{R}^ m$ for the contracting, expanding and isometric subspaces of the matrix $M_f$ in \eqref{eq:companion}. Then there exist unique points $y ^{\Delta _+}\in (W ^{(s)}+\mathbf{e}^{(1)})\cap (W ^{(u)}+W ^{(0)})$ and $y ^{\Delta _-}\in (W ^{(s)}+W ^{(0)}+\mathbf{e}^{(1)})\cap W ^{(u)}$, and $x ^{\Delta _\pm}=\pi (y ^{\Delta _\pm})$ (cf. \eqref{eq:homoclinic}--\eqref{eq:decay}).

We return to the general nonhyperbolic setting and put
\begin{equation}
\label{eq:l*}
\ell ^*(\mathbb{Z},\mathbb{R})=\biggl\{ w=(w_n) \in \mathbb{R}^ \mathbb{Z}: \sup_{n \in \mathbb{Z}} \frac{|w_n|}{|n|+1}<\infty \biggr\} \supset \ell ^ \infty (\mathbb{Z},\mathbb{R}),
\end{equation}
denote by $\ell ^*(\mathbb{Z},\mathbb{Z})$ the group of integer sequences in $\ell ^*(\mathbb{Z},\mathbb{R})$, and furnish these spaces with the topology of coordinate-wise convergence. We extend the maps $\bar{\sigma }$, $f(\bar{\sigma })$ and $\rho $ in \eqref{eq:sigmabar}--\eqref{eq:rho} to group homomorphisms
\begin{gather*}
\bar{\sigma }^*\colon \ell ^*(\mathbb{Z},\mathbb{R})\longrightarrow \ell ^*(\mathbb{Z},\mathbb{R}),\qquad f(\bar{\sigma }^*)\colon \ell ^*(\mathbb{Z},\mathbb{R})\longrightarrow \ell ^*(\mathbb{Z},\mathbb{R}),
\\
\rho ^*\colon \ell ^*(\mathbb{Z},\mathbb{R})\longrightarrow \mathbb{T}^ \mathbb{Z},
\end{gather*}
and set
\begin{equation}
\label{eq:Wf*}
\begin{aligned}
W_f ^*&=\{ w \in \ell ^*(\mathbb{Z},\mathbb{R}):f(\bar{\sigma }^*)(w)\in \ell ^*(\mathbb{Z},\mathbb{Z})\}
\\
&=\{ w \in \ell ^*(\mathbb{Z},\mathbb{R}):\rho ^*(w)\in X \}.
\end{aligned}
\end{equation}
Then $W_f ^*\subset \ell ^*(\mathbb{Z},\mathbb{R})$ is a closed, $\bar{\sigma }^*$-invariant subgroup. Since
$$
\ker f(\bar{\sigma }^*)=\{ w \in \ell ^*(\mathbb{Z},\mathbb{R}):f(\bar{\sigma }^*)(w)=0 \} \subset \ell ^ \infty (\mathbb{Z},\mathbb{R})
$$
(cf. \eqref{eq:womega}), we obtain that
$$
\ker f(\bar{\sigma }^*)=\ker f(\bar{\sigma })=W_f ^{(0)} \subset W_f.
$$

We define continuous group homomorphisms $\bar{\xi }^*\colon \ell ^ \infty (\mathbb{Z},\mathbb{Z})\longrightarrow W ^*_f$ and $\xi ^*\colon \ell ^ \infty (\mathbb{Z},\mathbb{Z})\longrightarrow X$ by setting
\begin{equation}
\label{eq:xi*}
\begin{gathered}
\bar{\xi }^*(v)=\sum_{n\ge0}v_n \bar{\sigma }^{-n}(w ^{\Delta _-})+\sum_{n<0}v_n \bar{\sigma }^{-n}(w ^{\Delta _+}),
\\
\xi ^*(v)=\rho ^* \circ \bar{\xi }^*(v)
\end{gathered}
\end{equation}
for every $v=(v_n)\in \ell ^ \infty (\mathbb{Z},\mathbb{Z})$. Since the coordinates $w_n ^{\Delta _+}$ and $w_{-n}^{\Delta _-}$ decay exponentially as $n \to \infty $ by \eqref{eq:decay}, each coordinate of $\bar{\xi }^*(v)$ in \eqref{eq:xi*} converges and
\begin{equation}
\label{eq:xi*2}
\bar{\xi }^*(\ell ^ \infty (\mathbb{Z},\mathbb{Z}))\subset W ^*_f.
\end{equation}
According to \eqref{eq:homoclinic2},
\begin{equation}
\label{eq:inverse2}
\begin{gathered}
f(\bar{\sigma }^*)\circ \bar{\xi }^*(v)=v,
\\
\bar{\xi }^*\circ f(\bar{\sigma }^*)(w)-w \in W_f ^{(0)}
\end{gathered}
\end{equation}
for every $v \in \ell ^ \infty (\mathbb{Z},\mathbb{Z})$ and $w \in W_f ^*$ (for the second equation in \eqref{eq:inverse2} we note that the maps $\bar{\sigma }^*$, $f(\bar{\sigma }^*)$ and $\bar{\xi }^*$ can be extended to the set of sequences with polynomial growth in $\mathbb{R}^ \mathbb{Z}$ and $\mathbb{Z}^ \mathbb{Z}$, respectively, where they still satisfy the first equation in \eqref{eq:inverse2};
%%%CHANGE
we note that the second equation also extends to such sequences). From \eqref{eq:Wf*} and \eqref{eq:inverse2} we see that
\begin{equation}
\label{eq:inverse3}
\begin{gathered}
\bar{\xi }^*\circ f(\bar{\sigma })(W_f)\subset W_f,
\\
V_f=f(\bar{\sigma })(W_f)=\{ v \in \ell ^ \infty (\mathbb{Z},\mathbb{Z}):\bar{\xi }^*(v)\in \ell ^ \infty (\mathbb{Z},\mathbb{R})\}.
\end{gathered}
\end{equation}

The map $\bar{\xi }^*\colon \ell ^ \infty (\mathbb{Z},\mathbb{Z})\longrightarrow \ell ^*(\mathbb{Z},\mathbb{R})$ can obviously not be expected to be shift-equivariant. Indeed,
\begin{equation}
\label{eq:commutation}
\begin{aligned}
\mathsf{d}(n,v)&=(\bar{\sigma }^*)^ n \circ \bar{\xi }^*(v)-\bar{\xi }^*\circ (\bar{\sigma }^*)^ n(v)
\\
&=
\begin{cases}
\sum_{j=0}^{n-1}v_j \bar{\sigma }^{n-j}w ^{\Delta _0}&\textup{if}\enspace n>0,
\\
0&\textup{if}\enspace n=0,
\\
-\sum_{j=1}^{n}v_{-j}\bar{\sigma }^{j-n}w ^{\Delta _0}&\textup{if}\enspace n<0.
\end{cases}
\end{aligned}
\end{equation}
for every $n \in \mathbb{Z}$ and $v \in \ell ^ \infty (\mathbb{Z},\mathbb{Z})$, and the resulting map
\begin{equation}
\label{eq:d}
\mathsf{d}\colon \mathbb{Z}\times \ell ^ \infty (\mathbb{Z},\mathbb{Z})\longrightarrow W_f ^{(0)}
\end{equation}
satisfies the cocycle equation
\begin{equation}
\label{eq:d-cocycle}
\mathsf{d}(m,\bar{\sigma }^ nv)+\bar{\sigma }^ m \mathsf{d}(n,v)=\mathsf{d}(m+n,v)
\end{equation}
for every $m,n \in \mathbb{Z}$ and $v \in \ell ^ \infty (\mathbb{Z},\mathbb{Z})$.

From the first formula in \eqref{eq:xi*} it is clear there exists a constant $c'>0$ with $|\bar{\xi }^*(v)|_n\le c'\cdot \| v \|_\infty $ for every $v \in \ell ^ \infty (\mathbb{Z},\mathbb{Z})$ and $n=0,\dots ,m-1$, where $m$ is the degree of $f$. Hence
$$
|\bar{\xi }^*\circ f(\bar{\sigma })(w)|_n\le c' \| f(\bar{\sigma })(w)\|_\infty \le c' \| f \|_1 \cdot \| w \|_\infty
$$
for every $w \in \ell ^ \infty (\mathbb{Z},\mathbb{R})$ and $n=0,\dots ,m-1$. Since there exists a constant $c''>0$ with
\begin{equation}
\label{eq:bound2}
\| w \|_\infty \le c''\cdot \max\,\{|w_0|,\dots ,|w_{m-1}|\}
\end{equation}
for every $w \in W_f ^{(0)} $ by \eqref{eq:womega}, we can find a constant $c>0$ with
\begin{equation}
\label{eq:bound}
\| \bar{\xi }^*\circ f(\bar{\sigma })(w)\|_\infty \le c \cdot \| f(\bar{\sigma })(w)\|_\infty \le c \cdot \| f \|_1 \cdot \| w \|_\infty
\end{equation}
for every $w \in \ell ^ \infty (\mathbb{Z},\mathbb{R})$.

Equation \eqref{eq:commutation} shows that the map $\xi ^*\colon \ell ^ \infty (\mathbb{Z},\mathbb{Z})\longrightarrow X$ is \emph{equivariant modulo $X ^{(0)}$}, and our next result implies that $\xi ^*$ is also \emph{surjective modulo $X ^{(0)}$}.

\begin{prop}
\label{p:2}
There exists a closed, bounded, shift-invariant subset $Y \subset V_f$ with $\xi ^*(Y)+X ^{(0)}=X$.
\end{prop}

\begin{proof}
We recall the notation
$B_r(\ell ^ \infty (\mathbb{Z},\mathbb{R})) = \left\{ w \in \ell ^ \infty (\ZZ, \RR ): \| w \| _ \infty \leq r \right\}$
(cf.~\eqref{eq:Br}) and set
\begin{equation}
\label{eq:BrWf}
B_r(W_f)=W_f \cap B_r(\ell ^ \infty (\mathbb{Z},\mathbb{R})).
\end{equation}
Then $B_r(W_f)$ is a closed and bounded --- and hence compact --- shift-invariant subset of $\ell ^ \infty (\mathbb{Z},\mathbb{R})$, and
\begin{equation}
\label{eq:Yr}
Y_r=f(\bar{\sigma })(B_r(W_f))
\end{equation}
is a compact shift-invariant subset of $V_f$.

For $r\ge1/2$, $\rho (B_r(W_f))=X$, and \eqref{eq:inverse2} guarantees that $\bar{\xi }^*\circ f(\bar{\sigma })(w)-w \in W_f ^{(0)} $ for every $w \in B_r(W_f)$. Hence there exists, for every $x \in X$, an element $w \in B_r(W_f)$ with $\rho (w)=x$, and the element $v=f(\bar{\sigma })(w)\in Y_r$ satisfies that $\bar{\xi }^*(v)-w \in W_f ^{(0)} $ and $\xi ^*(v)-x \in X ^{(0)}$. This proves (2).
\end{proof}

Proposition \ref{p:2} suggests the following definition.

\begin{defi}
\label{d:pseudocover}
A closed, bounded, shift-invariant subset $V \subset \ell ^ \infty (\mathbb{Z},\mathbb{Z})$ is a \emph{pseudo-cover} of $X$ if $\xi ^*(V)+X ^{(0)}=X$.
\end{defi}

The last part of this section is devoted to the question whether --- and to what extent --- the non-equivariance of $\xi ^*$ can be `corrected'. We start by showing that there is no continuous, equivariant and surjective map $\phi $ from $\ell ^ \infty (\mathbb{Z},\mathbb{Z})$ (or from any shift of finite type $Y \subset \ell ^ \infty (\mathbb{Z},\mathbb{Z})$) to $X$.

\begin{defi}
\label{d:homoclinic2}
Let $T$ be a homeomorphism of a compact metrizable space $Y$, and let $\delta $ be a metric on $Y$. Two points $x,y \in Y$ are \emph{homoclinic} if $\lim_{|n|\to \infty }\delta (T ^ nx,T ^ ny)=0$. The \emph{homoclinic equivalence relation} $\boldsymbol{\Delta }_T(Y)$ is defined as
$$
\boldsymbol{\Delta }_T(Y)=\{(x,y)\in Y ^ 2:x\enspace \textup{and}\enspace y\enspace \textup{are homoclinic}\}.
$$
For every $x \in Y$ we denote by
$$
\boldsymbol{\Delta }_T(x)=\{ y \in Y:(x,y)\in \boldsymbol{\Delta }_T(Y)\}
$$
the \emph{homoclinic equivalence class} of $x$. The homoclinic relation $\boldsymbol{\Delta }_T(Y)$ is \emph{topologically transitive} if $\boldsymbol{\Delta }_T(y)$ is dense in $Y$ for \emph{some} $y \in Y$, and \emph{minimal} if $\boldsymbol{\Delta }_T(y)$ is dense in $Y$ for \emph{every} $y \in Y$.

Note that all these definitions are independent of the specific choice of the metric $\delta $.
\end{defi}

\begin{prop}
\label{p:xi*}
Let $\alpha $ be an irreducible, ergodic and nonexpansive automorphism of a compact connected abelian group $X$, and let $T$ be a homeomorphism of a compact metrizable space $Y$ whose homoclinic relation $\boldsymbol{\Delta }_Y(T)$ is topologically transitive. If $\phi \colon Y \longrightarrow X$ a continuous equivariant map then $\phi (Y)$ consists of a single fixed point $\bar{x}$ of $\alpha $ in $X$.
\end{prop}

\begin{proof}
For any pair $y,y'$ of homoclinic points in $Y$, the points $\phi (y)$ and $\phi (y')$ are homoclinic in $X$, and hence $\phi (y)-\phi (y')=0$ by Corollary \ref{c:1}. If $\boldsymbol{\Delta }_T(y)$ is dense in $Y$ for some $y \in Y$ then the continuity of $\phi $ implies that $\phi (Y)$ is a single point which must be fixed under $\alpha $.
\end{proof}

\begin{coro}
\label{c:xi*}
Let $f \in R_1$ be an irreducible nonhyperbolic polynomial which is not cyclotomic, $\alpha =\alpha _{R_1/(f)}$ the ergodic and nonexpansive automorphism of the compact connected abelian group $X=X_{R_1/(f)}$ defined in \eqref{eq:principal}--\eqref{eq:alpha2}, and let $\bar{\sigma }$ be the shift \eqref{eq:sigmabar} on $\ell ^ \infty (\mathbb{Z},\mathbb{Z})$. If $\phi \colon \ell ^ \infty (\mathbb{Z},\mathbb{Z})\longrightarrow X$ is a continuous equivariant map, then $\phi (\ell ^ \infty (\mathbb{Z},\mathbb{Z}))$ consist of a single fixed point of $\alpha $.
\end{coro}

\begin{proof}
For every positive integer $r$ we set $B_r=\{ v \in \ell ^ \infty (\mathbb{Z},\mathbb{Z}):\| v \|_\infty \le r \}=\{-r,\dots ,r \} ^ \mathbb{Z}$. Then the restriction $T=\bar{\sigma }|_{B_r}$ of $\bar{\sigma }$ to $B_r$ has a topologically transitive homoclinic equivalence relation, and Proposition \ref{p:xi*} implies that $\phi (B_r)$ consists of a single point. Since this is true for every $r\ge1$, $\phi (\ell ^ \infty (\mathbb{Z},\mathbb{Z}))$ consist of a single fixed point of $\alpha $.
\end{proof}

\begin{defi}
\label{d:sft}
Let $A$ be a countably infinite set, $P=(P(a,a'),\,a,a'\in A))$ a transition matrix with entries in $\{ 0,1 \}$, and $X_P=\{ x=(x_n)\in A ^ \mathbb{Z}:P(x_n,x_{n+1})=1\enspace \textup{for every}\enspace n \in \mathbb{Z}\}$ the associated shift of finite type.

The shift $X_P$ and the transition matrix $P$ are \emph{irreducible} if there exist, for every $a,a'\in A$, an $n\ge1$ and elements $a_1=a,a_2,\dots ,a_n=a'$ in $A$ with $P(a_i,a_{i+1})=1$ for $i=1,\dots ,n-1$. If $X_P$ is irreducible then the \emph{period} $p(X_P)$ is the highest common factor of the set of integers $n\ge1$ for which there exist elements elements $a_1,a_2,\dots ,a_n=a_1$ in $A$ with $P(a_i,a_{i+1})=1$ for $i=1,\dots ,n-1$.

The shift $X_P$ and the matrix $P$ are \emph{mixing} if they are irreducible with period $1$.
\end{defi}

For the following corollaries we assume that $\alpha $ is an irreducible, ergodic and nonexpansive automorphism of a compact connected abelian group $X$.

\begin{coro}
\label{c:xi*1}
Let $Y$ be a mixing shift of finite type with finite or countably infinite alphabet. Then every continuous equivariant map $\phi \colon Y \longrightarrow X$ sends $Y$ to a single point.
\end{coro}

\begin{proof}
If $T$ is the shift on $Y$, then $\boldsymbol{\Delta }_Y(T)$ is minimal, and our claim follows from Proposition \ref{p:xi*}.
\end{proof}

\begin{coro}
\label{c:xi*2}
Let $Y$ be an irreducible shift of finite type with finite or countably infinite alphabet. Then every continuous equivariant map $\phi \colon Y \longrightarrow X$ sends $Y$ to a finite set.
\end{coro}

\begin{proof}
If $p(Y)$ is the period of $Y$ then the shift $T$ on $Y$ has the property that there exists a closed subset $Y_0 \subset Y$ such that the sets $T ^ kY_0,\,k=0,\dots ,n-1$ are disjoint, $T ^ nY_0=Y_0$, $\bigcup_{k=0}^{n-1}T ^ kY_0=Y$, and the restriction of $T ^ n$ to $Y_0$ is mixing. Hence $\phi (Y)$ consists of a single periodic orbit for any continuous equivariant map $\phi \colon Y \longrightarrow X$.
\end{proof}

\begin{coro}
\label{c:xi*4}
Let $Y$ be a topologically transitive sofic shift with finite alphabet. Then every continuous equivariant map $\phi \colon Y \longrightarrow X$ sends $Y$ to a finite set.
\end{coro}

\begin{proof}
Since $Y$ is a continuous equivariant image of a topologically transitive shift of finite type with finite alphabet, the result follows from Corollary \ref{c:xi*2}.
\end{proof}

\begin{rema}
\label{r:xi*}
Some non-sofic shift-spaces have a topologically transitive homoclinic equivalence relation. For example, in Proposition~\ref{p:beta} (3) we show that if $\beta >1$ is a real number, and if $V_\beta \subset \{ 0,\dots ,\lceil \beta -1 \rceil \} ^ \mathbb{Z}$ is the two-sided beta-shift space defined in \eqref{eq:Vbeta}, then the homoclinic equivalence relation $\boldsymbol{\Delta }_{\bar \sigma}(V_\beta )$ of the beta-shift $\sigma _\beta $ is topologically transitive, although $V_\beta $ is in general not sofic.
\end{rema}

\emph{The Corollaries \ref{c:xi*1}--\ref{c:xi*4} and Remark \ref{r:xi*} imply that $\alpha $ cannot have Markov \textup{(}or sofic\textup{)} partitions or covers in any reasonable sense, and that there are no nontrivial continuous equivariant maps from beta-shifts to $X$.}

\medskip In order to understand to what extent the non-equivariance of $\xi ^*$ can be `corrected' if we are allowed to drop continuity we set
\begin{equation}
\label{eq:L}
\tilde{Y}=\ell ^ \infty (\mathbb{Z},\mathbb{Z})\times W_f ^{(0)}\cong \ell ^ \infty (\mathbb{Z},\mathbb{Z})\times X ^{(0)}
\end{equation}
(where we are using the fact that the restriction $\rho |_{W_f ^{(0)}}\colon W_f ^{(0)}\longrightarrow X ^{(0)}$ is a bijection) and consider the continuous surjective maps $\tilde{\sigma }\colon \tilde{Y}\longrightarrow \tilde{Y}$ and $\tilde{\xi }^*\colon \tilde{Y}\longrightarrow W_f ^*$, defined by
\begin{equation}
\label{eq:tildemaps}
\begin{gathered}
\tilde{\sigma }(v,w)=(\bar{\sigma }v,\bar{\sigma }w+\mathsf{d}(1,v)),
\\
\tilde{\xi }^*(v,w)=\bar{\xi }^*(v)+w
\end{gathered}
\end{equation}
for every $(v,w)\in \tilde{Y}=\ell ^ \infty (\mathbb{Z},\mathbb{Z})\times W_f ^{(0)}$. The map $\tilde{\sigma }$ is obviously a homeomorphism, and
\begin{equation}
\label{eq:equivariance2}
\tilde{\xi }^*\circ \tilde{\sigma }=\bar{\sigma }^*\circ \tilde{\xi }^*.
\end{equation}
Finally we write $\tilde{\pi }\colon \tilde{Y}\longrightarrow \ell ^ \infty (\mathbb{Z},\mathbb{Z})$ for the first coordinate projection.

\begin{theo}
\label{t:3}
Let $f \in R_1$ be an irreducible nonhyperbolic polynomial which is not cyclotomic, $\alpha =\alpha _{R_1/(f)}$ the ergodic and nonexpansive automorphism of the compact connected abelian group $X=X_{R_1/(f)}$ defined in \eqref{eq:principal}--\eqref{eq:alpha2}, and let $\tilde{\sigma }\colon \tilde{Y}\longrightarrow \tilde{Y}$ be defined by \eqref{eq:L}--\eqref{eq:tildemaps}. For every $\bar{\sigma }$-invariant probability measure $\nu $ on $\ell ^ \infty (\mathbb{Z},\mathbb{Z})$ the following conditions are equivalent.

\begin{enumerate}
\item
There exists a $\tilde{\sigma }$-invariant probability measure $\tilde{\nu }$ on $\tilde{Y}$ with $\tilde{\pi }_*\tilde{\nu }=\nu $;
\item
For every $\varepsilon >0$ there exists a compact subset $C_\varepsilon \subset W_f ^{(0)}$ with
\begin{equation}
\label{eq:w-bounded}
\nu (\{ v \in \ell ^ \infty (\mathbb{Z},\mathbb{Z})\}:\mathsf{d}(k,v)\in C_\varepsilon \})>1-\varepsilon \enspace \textup{for every}\enspace k \in \mathbb{Z};
\end{equation}
\item
There exists a Borel map $\mathsf{b}\colon \ell ^ \infty (\mathbb{Z},\mathbb{Z})\longrightarrow W_f ^{(0)} $ with
\begin{equation}
\label{eq:b1}
\mathsf{d}(1,v)=\mathsf{b}(\bar{\sigma }v)-\bar{\sigma }\mathsf{b}(v)\enspace \textup{for}\enspace \nu \textsl{-a.e.}\;v \in \ell ^ \infty (\mathbb{Z},\mathbb{Z}).
\end{equation}
\end{enumerate}
If $\nu $ satisfies these equivalent conditions, then the Borel map $\xi ^*_\mathsf{b}\colon \ell ^ \infty (\mathbb{Z},\mathbb{Z})\linebreak[0]\longrightarrow X$, defined by
\begin{equation}
\label{eq:xi*b}
\xi ^*_\mathsf{b}(v)=\xi (v)+\rho ^*\circ \mathsf{b}(v)
\end{equation}
for every $v \in \ell ^ \infty (\mathbb{Z},\mathbb{Z})$, has the property that
\begin{equation}
\label{eq:xi*b2}
\begin{gathered}
\xi ^*_\mathsf{b}(v)-\xi (v)\in X ^{(0)}\enspace \textup{for every}\enspace v \in \ell ^ \infty (\mathbb{Z},\mathbb{Z}),
\\
\xi ^*_\mathsf{b} \circ \bar{\sigma }=\alpha \circ \xi ^*_\mathsf{b}\enspace \nu \textsl{-a.e.},
\end{gathered}
\end{equation}
and the probability measure
\begin{equation}
\label{eq:nuxi}
\mu =(\xi ^*_\mathsf{b})_*\nu
\end{equation}
on $X$ is $\alpha $-invariant.
\end{theo}

Motivated by Theorem \ref{t:3} we adopt the following terminology.

\begin{defi}
\label{d:bounded}
A shift-invariant probability measure $\nu $ on $\ell ^ \infty (\mathbb{Z},\mathbb{Z})$ is \emph{weakly $\mathsf{d}$-bounded} if it satisfies the three equivalent conditions of Theorem \ref{t:3}. The probability measure $\nu $ is \emph{$\mathsf{d}$-bounded} if there exists a compact subset $C \subset W_f ^{(0)}$ such that
\begin{equation}
\label{eq:bounded}
\nu (\{ v \in \ell ^ \infty (\mathbb{Z},\mathbb{Z}):\mathsf{d}(k,v)\in C\enspace \textup{for every}\enspace k \in \mathbb{Z}\})=1.
\end{equation}
\end{defi}

\begin{proof}[Proof of Theorem \ref{t:3}]
Suppose that $\tilde{\nu }$ is a $\tilde{\sigma }$-invariant probability measure on $\tilde{Y}$. We set $\nu =\tilde{\pi }_*\tilde{\nu }$, fix $\varepsilon >0$ and choose $K_\varepsilon >0$ with $\tilde{\nu }(\ell ^ \infty (\mathbb{Z},\mathbb{Z})\times B(K_\varepsilon ))>1-\varepsilon /2$, where $B(K_\varepsilon )=\{ w \in W_f ^{(0)}:\| w \|_\infty <K_\varepsilon \}$. Since $\tilde{\nu }$ is $\tilde{\sigma }$-invariant, $\tilde{\nu }\bigl(\tilde{\sigma }^ k(\ell ^ \infty (\mathbb{Z},\mathbb{Z})\times B(K_\varepsilon ))\cap (\ell ^ \infty (\mathbb{Z},\mathbb{Z})\times B(K_\varepsilon ))\bigr)>1-\varepsilon $ for every $k \in \mathbb{Z}$, which implies that
$$
\nu (\{ v \in \ell ^ \infty (\mathbb{Z},\mathbb{Z}):\| \mathsf{d}(k,v)\|_\infty <2K_\varepsilon \})>1-\varepsilon \enspace \textup{for every}\enspace k \in \mathbb{Z}.
$$
Since $\varepsilon $ was arbitrary this shows that (1) $\Rightarrow $ (2).

In order to check that (2) $\Rightarrow $ (3) we choose an enumeration $\Omega _f ^{(0)} =\{ \omega _1,\dots ,\omega _{m_0}\}$ of $\Omega _f ^{(0)} $, write $\mathbf{W}_f ^{(0)} = \mathbb{C}\otimes_\mathbb{R}W_f ^{(0)} $ for the complexification of $W_f ^{(0)} $ and use the basis $\{ w(\omega _i):i=1,\dots ,{m_0}\}$ in \eqref{eq:womega} to identify $\mathbf{W}_f ^{(0)}$ with $\mathbb{C}^{m_0}$. Let
\begin{equation}
\label{eq:Sm0}
\mathbb{S}^{m_0}=\{(\gamma _1,\dots ,\gamma _{m_0})\in \mathbb{C}^{m_0}:|\gamma _i|=1\enspace \textup{for}\enspace i=1,\dots ,{m_0}\}
\end{equation}
and define, for every $\gamma =(\gamma _1,\dots ,\gamma _{m_0})\in \mathbb{S}^{m_0}$, a linear map $M_\gamma \colon \mathbb{C}^{m_0}\longrightarrow \mathbb{C}^{m_0}$ by setting
\begin{equation}
\label{eq:Mgamma}
M_\gamma \mathbf{z}=(\gamma _1z_1,\dots ,\gamma _{m_0}z_{m_0})
\end{equation}
for every $\mathbf{z}=(z_1,\dots, z_{m_0})\in \mathbb{C}^{m_0}$. We form the locally compact semi-direct product
$$
\mathbf{G}=\mathbb{C}^{m_0}\rtimes \mathbb{S}^{m_0}
$$
with group operation
$$
(\mathbf{z},\gamma )\cdot (\mathbf{z}',\gamma ')=(\mathbf{z}+M_\gamma \mathbf{z}',\gamma \gamma ')
$$
for every $\mathbf{z},\mathbf{z}'\in \mathbb{C}^{m_0}$ and $\gamma ,\gamma '\in \mathbb{S}^{m_0}$ and set
$$
\mathsf{d} ^*(k,v)=(\mathsf{d}(k,v),\boldsymbol{\omega }^ k)
$$
for every $k \in \mathbb{Z}$ and $v \in \ell ^ \infty (\mathbb{Z},\mathbb{Z})$, where
\begin{equation}
\label{eq:pmb}
\boldsymbol{\omega }^ k=(\omega _1 ^ k,\dots ,\omega _{m_0}^ k).
\end{equation}
By \eqref{eq:d-cocycle}, the resulting map $\mathsf{d} ^*\colon \mathbb{Z}\times \ell ^ \infty (\mathbb{Z},\mathbb{Z})\longrightarrow \mathbf{G}$ satisfies the cocycle equation
$$
\mathsf{d} ^*(k,\bar{\sigma }^ lv)\cdot \mathsf{d} ^*(l,v)=\mathsf{d} ^*(k+l,v)
$$
for every $k,l \in \mathbb{Z}$ and $v \in \ell ^ \infty (\mathbb{Z},\mathbb{Z})$, where we are using the fact that the shift $\bar{\sigma }$ on $\mathbf{W}_f ^{(0)} $ corresponds to $M_{\boldsymbol{\omega }}$ under our identification of $\mathbf{W}_f ^{(0)} $ with $\mathbb{C}^{m_0}$. Furthermore, the map $\mathsf{d} ^*(k,\cdot )\colon \ell ^ \infty (\mathbb{Z},\mathbb{Z})\longrightarrow \mathbf{G}$ is continuous for every $k \in \mathbb{Z}$.

If the probability measure $\nu $ satisfies (2), then the cocycle $\mathsf{d} ^*\colon \mathbb{Z}\times \ell ^ \infty (\mathbb{Z},\mathbb{Z})\linebreak[0]\longrightarrow \mathbf{G}$ is \emph{bounded} in the sense that there exists, for every $\varepsilon >0$, a compact subset $\mathbf{C}_\varepsilon \subset \mathbf{G}$ with
$$
\nu (\{ v \in \ell ^ \infty (\mathbb{Z},\mathbb{Z}):\mathsf{d} ^*(k,v)\notin \mathbf{C}_\varepsilon \})<\varepsilon
$$
for every $k \in \mathbb{Z}$, and \cite[Theorem 4.7]{S0} implies the existence of a Borel map $\mathsf{b}'\colon \ell ^ \infty (\mathbb{Z},\mathbb{Z})\longrightarrow \mathbf{G}$ and of a compact subgroup $\mathbf{K}\subset \mathbf{G}$ such that
\begin{equation}
\label{eq:K1}
\mathsf{b}'(\bar{\sigma }v)^{-1}\cdot \mathsf{d} ^*(1,v)\cdot \mathsf{b}'(v)\in \mathbf{K}
\end{equation}
for $\nu \textsl{-a.e.}\;v \in \ell ^ \infty (\mathbb{Z},\mathbb{Z})$.

Every compact subgroup $\mathbf{K}\subset \mathbf{G}$ is of the form
$$
\mathbf{K}=\{(w(\gamma ),\gamma ):\gamma \in \Gamma _0 \}
$$
for some compact subgroup $\Gamma _0 \subset \mathbb{S}^{m_0}$ and some Borel map $w\colon \Gamma _0 \longrightarrow \mathbb{C}^{m_0}$ satisfying the cocycle equation
$$
w(\gamma \gamma ')=w(\gamma )+M_{\gamma }w(\gamma ')
$$
As $\Gamma _0$ is compact, this cocycle is a coboundary, i.e. there exists a $\mathbf{t}\in \mathbb{C}^{m_0}$ with
$$
w(\gamma )=M_\gamma \mathbf{t}-\mathbf{t}
$$
for every $\gamma \in \Gamma _0$, and
\begin{equation}
\label{eq:K2}
\mathbf{K}=\{(M_\gamma \mathbf{t}-\mathbf{t},\gamma ):\gamma \in \Gamma _0 \}.
\end{equation}

We write the map $\mathsf{b}'$ in \eqref{eq:K1} as $\mathsf{b}'=(b_1,b_2)$ with $b_1\colon \ell ^ \infty (\mathbb{Z},\mathbb{Z})\longrightarrow \mathbb{C}^{m_0}$ and $b_2\colon \ell ^ \infty (\mathbb{Z},\mathbb{Z})\longrightarrow \mathbb{S}^{m_0}$. According to \eqref{eq:K1}--\eqref{eq:K2},
\begin{align*}
~&(b_1(\bar{\sigma }v),b_2(\bar{\sigma }v))^{-1}\cdot (\mathsf{d}(1,v),\boldsymbol{\omega })\cdot (b_1(v),b_2(v))
\\
&\enspace =(-M_{b_2(\bar{\sigma }v)}^{-1}b_1(\bar{\sigma }v)+M_{b_2(\bar{\sigma }v)}^{-1}\mathsf{d}(1,v)+M_{b_2(\bar{\sigma }v)^{-1} \boldsymbol{\omega }}b_1(v),b_2(\bar{\sigma }v)^{-1}\boldsymbol{\omega }b_2(v))
\\
&\enspace =(M_{b_2(\bar{\sigma }v)^{-1}\boldsymbol{\omega }b_2(v))}\mathbf{t}-\mathbf{t},b_2(\bar{\sigma }v)^{-1}\boldsymbol{\omega }b_2(v)))
\end{align*}
with $b_2(\bar{\sigma }v)^{-1}\boldsymbol{\omega }b_2(v))\in \Gamma _0$, and hence
\begin{align*}
\mathsf{d}(1,v)&=M_{\boldsymbol{\omega }}(M_{b_2(v)}\mathbf{t}-b_1(v)) - (M_{b_2(\bar{\sigma }v)}\mathbf{t}-b_1(\bar{\sigma }v))
\\
&=\bar{\sigma }(M_{b_2(v)}\mathbf{t}-b_1(v)) - (M_{b_2(\bar{\sigma }v)}\mathbf{t}-b_1(\bar{\sigma }v))
\end{align*}
for $\nu \textsl{-a.e.}\;v \in \ell ^ \infty (\mathbb{Z},\mathbb{Z})$. We set $b(v)=b_1(v)-M_{b_2(v)}\mathbf{t}$ for every $v \in \ell ^ \infty (\mathbb{Z},\mathbb{Z})$, view $b$ as a map from $\ell ^ \infty (\mathbb{Z},\mathbb{Z})$ to $\mathbf{W}_f ^{(0)} \supset W_f ^{(0)} $, and obtain that
$$
\mathsf{d}(1,v)=b(\bar{{\sigma }}v)-\bar{\sigma }b(v)
$$
for $\nu\textsl{-a.e.}\;v \in \ell ^ \infty (\mathbb{Z},\mathbb{Z})$. If
$$
\mathsf{b}(v)=(b(v)+\overline{b(v)})/2,
$$
where the bar denotes complex conjugation in $\mathbf{W}_f ^{(0)} =\mathbb{C} \otimes_\mathbb{R}W_f ^{(0)} $, then the resulting map $\mathsf{b}\colon \ell ^ \infty (\mathbb{Z},\mathbb{Z})\longrightarrow W_f ^{(0)} $ satisfies \eqref{eq:b1}, and \eqref{eq:xi*b2} follows from \eqref{eq:b1} and \eqref{eq:xi*b}. This completes the proof that (2) $\Rightarrow $ (3).

Finally, if (3) is satisfied, then there exists a unique, and obviously $\tilde{\sigma }$-invariant, probability measure $\tilde{\nu }$ on $\tilde{Y}$ with $\tilde{\pi }_*\tilde{\nu }=\nu $ and $\tilde{\nu }(\{(v,\mathsf{b}(v)):v \in \ell ^ \infty (\mathbb{Z},\mathbb{Z})\})=1$, which proves that (3) $\Rightarrow $ (1).

The final assertions \eqref{eq:xi*b2} and the $\alpha $-invariance of the probability measure $\mu $ in \eqref{eq:nuxi} are immediate consequences of \eqref{eq:b1}.
\end{proof}

We recall the following definition from \cite{LS1}.

\begin{defi}
\label{d:equivalent}
Let $f \in R_1$ be an irreducible nonhyperbolic polynomial which is not cyclotomic, and let $\alpha =\alpha _{R_1/(f)}$ be the ergodic and nonexpansive automorphism of the compact connected abelian group $X=X_{R_1/(f)}$ in \eqref{eq:principal}-\eqref{eq:alpha2}. Two $\alpha $-invariant probability measures $\mu_1,\mu_2$ on $X$ are \emph{centrally equivalent} if they have an invariant joining $\varrho $ (i.e. an $(\alpha \times \alpha)$-invariant measure $\varrho $ on $X \times X$ which projects to $\mu_1$ and $\mu_2$, respectively) so that, for $\varrho\textsl{-a.e.}\;(x,y)\in X \times X$, $x$ and $y$ lie on the same central leaf. In other words,
$$
x-y \in X ^{(0)}\enspace \textup{for}\enspace \varrho \textsl{-a.e.}\;(x,y)\in X \times X,
$$
where $X ^{(0)}\subset X$ is the central subgroup of $\alpha $ defined in \eqref{eq:X0}.
\end{defi}

It is not hard to show that any two centrally equivalent probability measures have the same entropy under $\alpha $. Since Lebesgue measure is the unique measure of maximum entropy for $\alpha$, it follows that the only measure centrally equivalent to Lebesgue measure is Lebesgue measure itself.

\begin{exam}
\label{e:measures}
If $\nu $ is a weakly $\mathsf{d}$-bounded $\bar{\sigma }$-invariant probability measure on $\ell ^ \infty (\mathbb{Z},\linebreak[0]\mathbb{Z})$, and if $\mathsf{b}',\mathsf{b}''\colon \ell ^ \infty (\mathbb{Z},\mathbb{Z})\longrightarrow W_f ^{(0)}$ are two maps satisfying \eqref{eq:b1}, then the $\alpha $-invariant probability measures $\mu =(\xi ^*_{\mathsf{b}'})_*\nu $ and $\mu '=(\xi ^*_{\mathsf{b}''})_*\nu $ are centrally equivalent, since $\xi ^*_{\mathsf{b}'}(x)-\xi ^*_{\mathsf{b}''}(x)\in X ^{(0)}$ for every $x \in X$.
\end{exam}

If the equation \eqref{eq:b1} has a measurable solution $\mathsf{b}$, then this solution is generally not unique. Given a weakly $\mathsf{d}$-bounded $\bar \sigma$-invariant probability measure we may thus try to choose $\mathsf{b}$ so that $(\xi ^{*} _ \mathsf{b}) _ {*} \nu$ is as simple as possible.

\begin{prop}
\label{p: choice of coboundary}
Let $\nu$ be a weakly $\mathsf{d}$-bounded $\bar \sigma$-invariant probability measure on $\ell ^ \infty (\ZZ, \ZZ)$, and let $\mathsf{b}$ be a solution of \eqref{eq:b1}. Assume that $(\xi ^{*} _ {\mathsf{b}}) _ {*} \nu$ is singular with respect to Lebesgue measure. Then there is a solution $\mathsf{b} '$ of \eqref{eq:b1} and an $\alpha$-invariant Borel set $ Z \subset X$ so that
\begin{enumerate}
\item
$Z$ intersect each coset of $X ^ {( 0 )}$ in at most one point,

\item
$(\phi ^{*} _ {\mathsf{b} '}) _ {*} \nu (Z) = 1$.
\end{enumerate}
\end{prop}

\begin{proof}
By \cite[Theorem 1.3.(1)] {LS1}, there exist a probability measure $\mu '$ on $X$ which is centrally equivalent to $ \mu = (\xi ^{*} _ {\mathsf{b}}) _ {*} \nu$ and a Borel set $Z \subset X$ (which we may as well assume to be $\alpha$-invariant) of full $\mu '$-measure which intersects each coset of $X ^ {( 0 )}$ in at most one point.

Since $\mu '$ and $\mu$ are centrally equivalent,
\begin{equation}
\label{equation: single point in intersection}
\bigl|{(\xi ^{*} _ {\mathsf{b}} (v) + X ^ {( 0 )}) \cap Z}\bigr|=1
\end{equation}
for $\nu\emph{-a.e.}\;v \in \ell ^ \infty (\mathbb{Z},\mathbb{Z})$. Define $\mathsf{b}'(v) \in W _ f ^ {( 0 )}$ by the requirement that
\begin{equation*}
\xi ^{*} (v) + \rho \circ \mathsf{b} ' (v)
\end{equation*}
is the single point in the set \eqref{equation: single point in intersection}. This function is certainly measurable, as can be verified easily by using the joining establishing the central equivalence of $\mu '$ and $\mu$. It also satisfies \eqref{eq:b1}:
since $Z$ is $\alpha$-invariant,
\begin{align*}
\left\{ \alpha \circ \xi ^{*} (v)+ \rho \circ \bar \sigma (\mathsf{b}'(v)) \right\} & =
\alpha (\xi ^{*} _ {\mathsf{b}} (v) + X ^ {( 0 )}) \cap Z
\\
& = (\xi ^{*} _ {\mathsf{b}} (\bar \sigma v) + X ^ {( 0 )}) \cap Z = \left\{ \xi ^{*} (\bar \sigma x) + \rho \circ \mathsf{b} (\bar \sigma x) \right\},
\end{align*}
and since $\rho$ is injective on $W ^ {( 0 )} _ f$,
\begin{equation*}
\mathsf{d} (1,v)= \bar \sigma ^* \circ \bar \xi ^{*} (v) - \bar \xi ^{*} \circ \bar \sigma (v) = \mathsf{b} ' (\bar \sigma v) - \bar \sigma \mathsf{b} '(v).
\end{equation*}
By construction, $\xi _ {\mathsf{b} '} (v) \in Z$ for every $v$ for which $\mathsf{b} '$ is well-defined (i.e. on a set of full $\nu $-measure).
\end{proof}

\begin{prop}
\label{p:measures}
Let $f \in R_1$ be an irreducible nonhyperbolic polynomial which is not cyclotomic, and let $\alpha =\alpha _{R_1/(f)}$ be the ergodic and nonexpansive automorphism of the compact connected abelian group $X=X_{R_1/(f)}$ in \eqref{eq:principal}-\eqref{eq:alpha2}. For every $\alpha $-invariant probability measure $\mu $ on $X$ there exists a $\mathsf{d}$-bounded $\bar{\sigma }$-invariant probability measure $\nu $ on $\ell ^ \infty (\mathbb{Z},\mathbb{Z})$ such that $\mu $ is centrally equivalent to the probability measure $(\xi ^*_\mathsf{b})_*\nu $ in Theorem \ref{t:3}.
\end{prop}

\begin{proof}
Let $\mu $ be an $\alpha $-invariant probability measure on $X$. We set $W=\{ w=(w_n)\in W_f:0\le w_n<1\enspace \textup{for every}\enspace n \in \mathbb{Z}\}$ (cf. \eqref{eq:Wf}), note that the restriction $\rho |_W$ of the equivariant map $\rho \colon W_f \longrightarrow X$ to $W$ is bijective, and conclude that there exists a unique $\bar{\sigma }$-invariant probability measure $\mu '$ on $W$ with $\rho _*\mu '=\mu $.

The $\bar{\sigma }$-invariant probability measure $\nu =f(\bar{\sigma })_*\mu '$ is supported on $Y_r \subset \ell ^ \infty (\mathbb{Z},\mathbb{Z})$ for some $r>0$, where $Y_r=f(\bar{\sigma })(B_r(W_f))$ as in the proof of Proposition \ref{p:2}. By \eqref{eq:bound}, $\bar{\xi }^*(Y_r)$ is a bounded subset of $\ell ^ \infty (\mathbb{Z},\mathbb{R})$.
This shows that the cocycle $\mathsf{d}$ in \eqref{eq:d} is uniformly bounded on $Y_r$ and so $\nu $ is $\mathsf{d}$-bounded (cf. \eqref{eq:Yr}).

Let $\mathsf{b}\colon \ell ^ \infty (\mathbb{Z},\mathbb{Z})\longrightarrow W_f ^{(0)}$ be a Borel map satisfying \eqref{eq:b1}, and let $\xi ^*_\mathsf{b}\colon \ell ^ \infty (\mathbb{Z},\mathbb{Z})\longrightarrow X$ be given by \eqref{eq:xi*b}. Since $x-\xi ^*_\mathsf{b} \circ f(\bar{\sigma })\circ (\rho |_W)^{-1}(x)\in X ^{(0)}$ for every $x \in X$, the $\alpha $-invariant probability measure $(\xi ^*_\mathsf{b})_*\nu $ is centrally equivalent to $\mu $.
\end{proof}

The discussion in this section shows that in the nonexpansive case we have to make a choice between continuity and equivariance: the map $\xi ^*\colon \ell ^ \infty (\mathbb{Z},\mathbb{Z})\linebreak[0]\longrightarrow X$ in \eqref{eq:xi*} is continuous, but not equivariant, and the maps $\xi _\mathsf{b}\colon \ell ^ \infty (\mathbb{Z},\mathbb{Z})\linebreak[0]\longrightarrow X$ in \eqref{eq:xi*b}, which are equivariant at least on some reasonably large sets, are generally not continuous. In neither case can we expect these maps to be surjective.

If $\nu $ is a weakly $\mathsf{d}$-bounded $\bar{\sigma }$-invariant probability measure on $\ell ^ \infty (\mathbb{Z},\mathbb{Z})$, then the Borel map $\xi _\mathsf{b}\colon \ell ^ \infty (\mathbb{Z},\mathbb{Z})\longrightarrow X$ in \eqref{eq:xi*b} is equivariant $\nu \emph{-a.e.}$ and the $\mu =(\xi ^*_\mathsf{b})_*\nu $ in \eqref{eq:nuxi} is therefore $\alpha $-invariant, but the entropy of $\mu $ will generally be lower than that of $\nu $.

In Proposition \ref{p:measures} we saw that we can obtain every $\alpha $-invariant probability measure on $X$ --- up to central equivalence --- from a $\mathsf{d}$-bounded shift-invariant probability measure on $\ell ^ \infty (\mathbb{Z},\mathbb{Z})$. However, all such measures are concentrated on the somewhat elusive set $V_f$, so that this result is of limited interest.

For this reason one would ideally like to find `nice' and `large' compact subshifts $V \subset \ell ^ \infty (\mathbb{Z},\mathbb{Z})$ (where \emph{nice} means something like a shift of finite type or a sofic shift, and \emph{large} means that the subshift should be a pseudo-cover of $X$ in the sense of Definition \ref{d:pseudocover}), such that the following conditions are satisfied:
\begin{enumerate}
\item
for every $v \in \ell ^ \infty (\mathbb{Z},\mathbb{Z})$, the intersection of $V$ with $v+\bigl(f(\bar{\sigma })(\ell ^ *(\mathbb{Z},\mathbb{Z}))\linebreak[0]\cap \ell ^ \infty (\mathbb{Z},\mathbb{Z})\bigr)$ is as small as possible,
\item
for every weakly $\mathsf{d}$-bounded shift-invariant probability measure $\nu $ on $V$, the probability measure $\mu =(\xi _\mathsf{b} ^*)_*\nu $ in \eqref{eq:nuxi} has the same entropy as $\nu $,
\item
every $\alpha $-invariant probability measure is centrally equivalent to a probability measure obtained in this manner.
\end{enumerate}

At this stage we have made only limited progress in this direction (cf. Section \ref{s:symbolic}, where we investigate the connection between two-sided beta-shift arising from a Salem number $\beta $ and the ergodic nonhyperbolic toral automorphism defined by the companion matrix if the minimal polynomial of $\beta $). One of the key difficulties one encounters in pursuing this program in any generality is the following: although the restriction to $\ell ^ \infty (\mathbb{Z},\mathbb{Z})$ of the map $f(\bar{\sigma })\colon W_f \longrightarrow \ell ^ \infty (\mathbb{Z},\mathbb{Z})$ is injective, the set
$$
\{ w \in \ell ^ *(\mathbb{Z},\mathbb{Z}):f(\bar{\sigma })(w)\in V \}
$$
need not be contained in $\ell ^ \infty (\mathbb{Z},\mathbb{Z})$ and the set
$$
\{ w \in \ell ^ \infty (\mathbb{Z},\mathbb{Z}):f(\bar{\sigma })(w)\in V \}
$$
may be unbounded, even if $V \subset \ell ^ \infty (\mathbb{Z},\mathbb{Z})$ is a bounded, shift-invariant set.

\begin{exam}
\label{e:unbounded}
For every $w \in W_f ^{(0)}$ and $n \in \mathbb{Z}$ we set $\zeta (w)_n=\lceil w_n \rceil$, where $\lceil t \rceil$ is again the smallest integer $\ge t$ for any $t \in \mathbb{R}$. The resulting map $\zeta \colon W_f ^{(0)}\longrightarrow \ell ^ \infty (\mathbb{Z},\mathbb{Z})$ has the property that the set $\zeta (W_f ^{(0)})$ is unbounded, but $\| f(\bar{\sigma })(v)\|_\infty <\| f \|_1$ for every $v \in \zeta (W_f ^{(0)})$.

In order to verify that $f(\bar{\sigma })$ maps some unbounded sequences in $\ell ^*(\mathbb{Z},\mathbb{Z})$ into $\ell ^ \infty (\mathbb{Z},\mathbb{Z})$ we choose $\theta \in \Omega _h ^{(0)}$ (cf. \eqref{eq:Omega}) and define, for every integer $j\ge0$, a point $\omega ^{(j)}=(\omega _n ^{(j)})\in \ell ^ \infty (\mathbb{Z},\mathbb{R})$ by setting
$$
\smash[b]{\omega _n ^{(j)}=
\begin{cases}
\theta ^ n+\theta ^{-n}&\textup{if}\enspace n\ge j,
\\
0&\textup{otherwise}.
\end{cases}}
$$
Then
$$
(f(\bar{\sigma })\omega ^{(j)})_n=0
$$
for $n<j-d$ and $n\ge j$, and $\| h(\bar{\sigma })\omega ^{(j)}\|_\infty \le 2 \cdot \| f \|_1$. For every $n \in \mathbb{Z}$ we put
$$
\tilde{\omega }=\sum_{j=0}^ \infty \omega ^{(3jd)}, \qquad w_n=\lceil \tilde{\omega }_n \rceil.
$$
The resulting point $w$ in $\ell ^*(\mathbb{Z},\mathbb{Z})$ is unbounded and satisfies that $\| f(\bar{\sigma })w \|_\infty \linebreak[0]\le 3 \cdot \| f \|_1$.
\end{exam}

\section{Beta-shifts and their properties}
\label{s: beta}

We fix a real number $\beta >1$ and consider the \emph{beta-transformation}
\label{discussion2}
$$
x \mapsto T_\beta x=\beta x \pmod 1
$$
from the closed unit interval $[0,1]$ to the half-open interval $I=[0,1)$ (cf. \cite{Parry} and \cite{Renyi}).

For every $x \in I$, the \emph{beta-expansion} $e_\beta (x)=(e_\beta (x)_n,\,n\ge1)$ of $x$ is defined by
$$
e_\beta (x)_n=\beta T_\beta ^{n-1}x-T_\beta ^ nx
$$
for every $n\ge1$. Note that $e_\beta (x)_n \in \{ 0,\dots ,\lceil \beta-1 \rceil \}$ for every $n\ge1$, where $\lceil t \rceil$ is the smallest integer $\ge t$ for any $t \in \mathbb{R}$, and that
\begin{equation}
\label{eq:expansion2}
x=\smash[b]{\sum_{n\ge1}e_\beta (x)_n \beta ^{-n}}
\end{equation}
for every $x \in I$.

We denote by $\prec$ the lexicographic order on the space $\ell ^ \infty _+$ of all bounded one-sided sequences $v=(v_n,\,n\ge1)$ of nonnegative integers and write $\bar{\sigma }_+$ for the one-sided shift $(\bar{\sigma }_+v)_n\linebreak[0]=v_{n+1}$ on $\ell ^ \infty _+$. The closed, $\bar{\sigma }_+$-invariant set
\begin{equation}
\label{eq:Vbeta+}
V_\beta ^+=\overline{\left\{ e _ \beta (x): x \in I \right\}}
\end{equation}
is called the \emph{beta-shift space} (where the bar denotes closure); it contains a unique lexicographically maximal element $e_\beta ^*$ with the following properties (cf. \cite{Parry}):
\begin{equation}
\label{eq:maximal}
\begin{gathered}
\bar{\sigma }_+^ k e_\beta ^*\preceq e_\beta ^*\enspace \textup{for every}\enspace k\ge0,
\\
\sum_{n\ge1}e_\beta ^*\beta ^{-n}=1, \qquad \bar{\sigma }_+^ n e_\beta ^*\ne 0\enspace \textup{for every}\enspace n\ge0,
\\
V_\beta ^+=\{ v \in \ell _+^ \infty :\bar{\sigma }_+^ nv\preceq e_\beta ^*\enspace \textup{for every}\enspace n\ge0 \}.
\end{gathered}
\end{equation}

Here we are interested in the two-sided beta-shift space. We write $v ^+=(v_1,v_2,\dots )$ for every $v=(v_n)\in \ell ^ \infty (\mathbb{Z},\mathbb{Z})$ and set
\begin{equation}
\label{eq:Vbeta}
\begin{aligned}
V_\beta &=\{ v \in \ell ^ \infty (\mathbb{Z},\mathbb{Z}):(\bar \sigma ^ nv)^+\in V_\beta ^+\enspace \textup{for every}\enspace n \in \mathbb{Z}\}
\\
&=\{ v \in \ell ^ \infty (\mathbb{Z},\mathbb{Z}):(\bar{\sigma }^ nv)^+\preceq e_\beta ^*\enspace \textup{for every}\enspace n \in \mathbb{Z}\}.
\end{aligned}
\end{equation}

For every $v \in \ell ^ \infty (\mathbb{Z},\mathbb{Z})$ with $v_{-n}=0$ for all sufficiently large $n\ge0$ we define the \emph{evaluation} $\eta _\beta (v)\in \mathbb{R}$ by
\begin{equation}
\label{eq:eta}
\eta _\beta (v)=\sum_{n \in \mathbb{Z}}v_n \beta ^{-n}.
\end{equation}
If we view $V_\beta ^+$ as the subset $\{ v \in V_\beta :v_n=0\enspace \textup{for}\enspace n\le0 \}$, then the evaluation defines a continuous, surjective, at most two-to-one map $\eta _\beta \colon V_\beta ^+\longrightarrow [0,1]$ with
\begin{equation}
\label{eq:unique0}
e_\beta (\eta _\beta (v))=v
\end{equation}
for all $v$ in the complement of a countable subset of $V_\beta ^+$ (cf. \cite{Parry} and \eqref{eq:expansion2}): the only possible exceptions to \eqref{eq:unique0} are points satisfying $(\bar{\sigma }^ kv)^+=e_\beta ^*$ for some $k>0$ (cf.\eqref{eq:maximal}).

The following elementary observations follow directly from \eqref{eq:maximal}--\eqref{eq:Vbeta}:

\begin{prop}
\label{p:beta}
Let $\beta >1$, and let $V_\beta \subset \ell ^ \infty (\mathbb{Z},\mathbb{Z})$ be the two-sided beta-shift space defined in \eqref{eq:Vbeta}.

\textup{(1)} If $v,w \in V_\beta $ satisfy that $w ^+\prec v ^+$ in the notation of \eqref{eq:maximal}--\eqref{eq:Vbeta}, then the point $v'$ with
$$
v_n'=
\begin{cases}
v_n&\textup{if}\enspace n\le0,
\\
w_n&\textup{if}\enspace n>0,
\end{cases}
$$
lies in $V_\beta $;

\textup{(2)} If $\beta $ is algebraic with minimal polynomial $f \in R_1$, and if $v,w \in V_\beta $ satisfy that $v_n=w_n$ for every $n<0$, $v_0>w_0$ and $w-v \in f(\bar{\sigma })(\ell ^ \infty (\mathbb{Z},\mathbb{Z}))$, then $v_0=w_0+1$, $v ^+=0$ and $w ^+=e_\beta ^*$. It follows that
$$
(v+f(\bar{\sigma })(\ell ^ 1(\mathbb{Z},\mathbb{Z})))\cap V_\beta =\{ v \}
$$
for every $v \in V_\beta $.

\textup{(3)} The homoclinic equivalence relation $\boldsymbol{\Delta } _ {\bar \sigma} (V _ \beta)$ \textup{(}cf. Definition~\ref{d:homoclinic2}\textup{)} is topologically transitive on $V_\beta$.
\end{prop}

\begin{proof}
In order to prove (1) we note that $(v_{k+1}',v_{k+2}',\dots )\prec (v_{k+1},v_{k+2},\dots )\linebreak[0]\preceq e_\beta ^*$ whenever $k<0$, and $(v_{k+1}',v_{k+2}',\dots )= (w_{k+1},w_{k+2},\dots )\preceq e_\beta ^*$ otherwise. According to \eqref{eq:Vbeta} this implies that $v'\in V_\beta $.

If $\beta $ is algebraic with minimal polynomial $f$, and if
$$
w \in (v+f(\bar{\sigma })(\ell ^ \infty (\mathbb{Z},\mathbb{Z})))\cap V_\beta
$$
and $v_n=w_n$ for all $n<0$, then
$$
\eta _\beta (w_0,w_1,\dots )=\eta _\beta (v_0,v_1,\dots ),
$$
and \eqref{eq:maximal} and \eqref{eq:unique0} imply (2).

In order to verify (3), we denote by $ \boldsymbol 0 \in V _ \beta$ the two-sided infinite sequence of zeros. For every $v \in V _ \beta$ and $n \in \NN$ the point $v '$ defined by
\begin{equation*}
v ' _ k =
\begin{cases}
0& \qquad \text{if $k < - n$}
\\
v _ n& \qquad \text{if $- n \leq k \leq n$}
\\
0& \qquad \text{if $n < k$}
\end{cases}
\end{equation*}
again lies in $V_\beta $, due to the lexicographic definition of the beta-shift in \eqref{eq:maximal}. It is also clearly homoclinic to $\boldsymbol 0$. This shows that the homoclinic equivalence class of $\boldsymbol 0$ is dense in $V_\beta $.
\end{proof}

Beta-shifts are in general not sofic (in fact, they are sofic if and only if the point $e_\beta ^*$ in \eqref{eq:maximal} is eventually periodic which implies, in turn, that $\beta $ is algebraic --- cf. e.g. \cite{Blanchard}). However, even if $V_\beta $ is not sofic, i.e. cannot be obtained by relabelling the letters of a shift of finite type with finite alphabet, it always has a nice description in terms of a certain shift of finite type $\Sigma _ \beta$ with a \emph{countable} alphabet. This infinite state shift of finite type has additional nice properties which make it a useful tool in the study of beta-shifts.

We now present the construction in \cite{Hofbauer, Takahashi} of this shift of finite type due to Hofbauer and Takahashi and its relation to the beta-shift (note that there is a gap in \cite{Takahashi}; see \cite{Hofbauer} for details).

For any pair of points $v, v' \in V _ \beta ^+$ with $v\prec v'$, let $[ v, v ' ] \subset V _ \beta ^+$ be the set of points which lie between $v$ and $v '$ in the lexicographic order, and let $\boldsymbol 0 _ +$ denote the one-sided sequence of zeros.

Let $A = \left\{ 0, \dots, \lceil \beta - 1 \rceil \right\}$. The alphabet (or state space) $\bar A$ of $ \Sigma _ \beta$ is given by $\bar A = \bar A ' \cup \bar A ''$ where
\begin{equation}
\label{eq:Abar}
\begin{aligned}
\bar A '& = \left\{ (a, [ \boldsymbol 0 _ +, e ^{*} _ \beta ]): a = 0, \dots, \lceil \beta - 2 \rceil \right\}
\\
\bar A ''& = \{ ((e ^{*} _ \beta) _ k, [ \boldsymbol 0 _ +, \bar \sigma _ + ^ k e ^{*} _ \beta ] ): k = 1, 2, \dots \}.
\end{aligned}
\end{equation}
Note that $\bar A$ is finite if and only if $e ^{*} _ \beta$ is eventually periodic.

The allowed transition in $\Sigma _\beta $ are defined as follows. Each state $\bar a \in \bar A '$ can be followed by any other state in $\bar A '$ as well as by the state $(\lceil \beta - 1 \rceil, [ \boldsymbol 0 _ +, \bar \sigma _ + e ^{*} _ \beta ] )$.
Each state $\bar a = ((e ^{*} _ \beta) _ k, [ \boldsymbol 0 _ +, \bar \sigma _ + ^ k e ^{*} _ \beta ] ) \in \bar A ''$ can be followed by either $\bar a ' = (a, [ \boldsymbol 0, e ^{*} _ \beta])$ for $a < (e ^{*} _ \beta) _ {k + 1}$, or by $((e ^{*} _ \beta) _ {k + 1}, [ \boldsymbol 0 _ +, \bar \sigma _ + ^ {k + 1} e ^{*} _ \beta ] ) $. We denote by
$P = (P (\bar a, \bar a '), \bar a, \bar a ' \in \bar A)$ the corresponding transition matrix, i.e.
$P (\bar a, \bar a ') = 1$ if and only if $ \bar a$ can be followed by $\bar a '$.

Let $\phi\colon \bar A \longrightarrow A$ be the projection onto the first coordinate, and let $\boldsymbol \phi$ be the corresponding map from ${\bar A} ^ \ZZ$ to $A ^ \ZZ$.  One can show quite easily that $\boldsymbol \phi (\Sigma _ \beta) \subset V _ \beta $.
In general, $\boldsymbol \phi |_ {\Sigma _ \beta}$ need not be surjective. What is true (see \cite{Hofbauer}) is that the complement $N$ of $\boldsymbol \phi (\Sigma _ \beta)$ is a shift invariant subset of $ V _ \beta$ with the property that any measure supported on it has zero entropy.

This construction is used in particular to show that $V _ \beta$ has a unique $\bar \sigma$-invariant measure $\mu _ \beta$ of maximal entropy with entropy $\log \beta$ (cf. \cite{Hofbauer}).

\begin{theo}
[\cite{Hofbauer}]
\label{t:Hofbauer}
The transition matrix $P=(P(a,a'),\,a,a'\in \bar{A})$ of $\Sigma _\beta $ has maximal eigenvalue $\beta $ and unique positive left and right eigenvectors $\mathbf{x}=(\mathbf{x}(a),\,a \in \bar{A})$, $\mathbf{y}=(\mathbf{y}(a),\,a \in \bar{A})$ with $\mathbf{x}P=\beta \mathbf{x}$, $P \mathbf{y}=\beta \mathbf{y}$ and $\sum_{a \in \bar{A}}\mathbf{x}(a)= \sum_{a \in \bar{A}}\mathbf{y}(a)=1$.

Let $\bar{\mu }_P$ be the Markov measure on $\Sigma _P$ defined by
\begin{equation}
\label{equation about Markov measure}
\bar{\mu }_P([a_{m_1},\dots ,a_{m_2}])=\frac{\beta ^{-(m_2-m_1)}\mathbf{y}(a_{m_2})}{\mathbf{y}(a_{m_1})}
\end{equation}
for every cylinder set
$$
[a_{m_1},\dots ,a_{m_2}]=\{ y \in \Sigma _P:y_n=a_n\enspace \textup{for}\enspace n=m_1,\dots ,\linebreak[0]m_2 \}.
$$
Then the restriction of $\boldsymbol{\phi }$ to $\Sigma _P$ is injective $\bar{\mu }_P\emph{-a.e.}$, and $\boldsymbol{\phi }_*\bar{\mu }_P=\mu _\beta $.
\end{theo}

The beta-shift is known to be sofic for Pisot numbers as well as for Salem numbers of degree four. For general Salem numbers $\beta $ it is not known whether $V_\beta$ has to be sofic (cf. \cite{Bertrand}, \cite{Boyd1}--\cite{Boyd3} and \cite{Sbeta}).

\section{The beta-shift and symbolic embeddings for Salem numbers}
\label{s:symbolic}

We start this section with a brief review of the case where $\beta$ is a Pisot number, and where the $\beta$-shift is a sofic model of the corresponding hyperbolic toral automorphism.

\begin{prop}[\cite{S3}]
\label{p:pisot}
Let $\beta >1$ be a Pisot number, $f \in R_1$ its minimal polynomial of degree $m$, say, and let $\alpha =\alpha _{R_1/(f)}$ be the expansive automorphism of the compact abelian group $X=X_{R_1/(f)}$ described in \eqref{eq:principal} \textup{(}if $\beta $ is a Pisot unit then $X \cong \mathbb{T}^ m$\textup{)}. Then the restriction of the equivariant map $\xi \colon \ell ^ \infty (\mathbb{Z},\mathbb{Z})\linebreak[0] \longrightarrow X$ in \eqref{eq:xi} to the two-sided beta-shift $V_\beta \subset \ell ^ \infty (\mathbb{Z},\mathbb{Z})$ is surjective and finite-to-one. In particular, if $\nu $ is a shift-invariant probability measure on $V_\beta $, then the measure $\mu =\xi _*\nu $ on $X$ is $\alpha $-invariant and has the same entropy as $\nu $. Furthermore, every $\alpha $-invariant probability $\mu $ on $X$ can be obtained in this manner.
\end{prop}
\label{one-to-one}

The restriction of $\xi $ to $V_\beta $ in Proposition \ref{p:pisot} is conjectured to be almost one-to-one, although this has only been proved in some examples (cf. \cite{S3}--\cite{Sid2}). For earlier special cases of Proposition \ref{p:pisot} we refer to \cite{SV}.

\smallskip
Proposition \ref{p:pisot} describes the close connection between the two-sided beta-shift of a Pisot unit $\beta >1$ and the toral automorphism defined by the companion matrix of the minimal polynomial of $\beta $. One of the principal motivations of this paper was the question whether there exists an analogous result for Salem numbers.

The following discussion shows that, although Proposition \ref{p:pisot} does not hold in this case, there \emph{does} exist a connection between two-sided beta-shifts of Salem numbers and the nonhyperbolic ergodic toral automorphisms defined by the companion matrices of their minimal polynomials. However, this connection is much more complicated and tenuous than in the Pisot case.

For the remainder of this section we restrict ourselves to Salem numbers, their minimal polynomials and their companion matrices. Assume therefore that $\beta $ is a Salem number with minimal polynomial $f \in R_1$ of (even) degree $m$, say, and let $\alpha =\alpha _{R_1/(f)}$ be the ergodic and nonexpansive automorphism of $X=X_{R_1/(f)}$ defined by \eqref{eq:principal}--\eqref{eq:alpha2}, which is algebraically conjugate to the companion matrix $M_f$ in \eqref{eq:companion}, acting on $\mathbb{T}^ m$. The corresponding two-sided beta-shift will be denoted by $V_\beta \subset \ell ^ \infty (\mathbb{Z},\mathbb{Z})$, and we write $\mu _\beta $ for the unique shift-invariant measure of maximal entropy on $V_\beta $.

Since the homoclinic equivalence relation $\boldsymbol{\Delta } _ {\bar \sigma} (V _ \beta)$ is topologically transitive on $V _ \beta$ by Proposition~\ref{p:beta} (3), Proposition~\ref{p:xi*} shows that a simple symbolic description as in Proposition~\ref{p:pisot} is not possible in this case. A partial analogue to Proposition~\ref{p:pisot} is presented in Theorem~\ref{t:weaklybounded} below.

%%%%%%%%%%%%%%%%%%%%%%%%%%%%%%%%%%

\begin{defi}
\label{definition: countable-to-one}
Let $Z_1, Z_2$ be standard Borel spaces and $\nu$ a probability measure on $Z _ 1$. A Borel map $g\colon Z _ 1 \longrightarrow Z _ 2$ is countable-to-one $\nu \textsl{-a.e.}$ if there are Borel sets $Z _ 1'\subset Z_1, Z _ 2 '\subset Z_2$ with $\nu (Z_1')=g _ {*} \nu(Z_2')=1$ so that $g ^{-1} (z) \cap Z _ 1 '$ is countable for every $z \in Z _ 2 '$.
\end{defi}

It is an easy exercise to see that entropy is preserved under almost everywhere countable-to-one factor maps.

\begin{theo}
\label{t:weaklybounded}
Let $\beta >1$ be a Salem number of degree $m$, say, $f \in R_1$ its minimal polynomial, and let $\alpha =\alpha _{R_1/(f)}$ be the ergodic and nonexpansive automorphism of $X=X_{R_1/(f)}\cong \mathbb{T}^ m$ defined in \eqref{eq:principal}--\eqref{eq:alpha2}.

Suppose that $\nu $ is a weakly $\mathsf{d}$-bounded $\bar{\sigma }$-invariant probability measure on the two-sided beta-shift $V_\beta $, and that $\xi _\mathsf{b} ^*\colon V_\beta \longrightarrow X$ is the $\nu \emph{-a.e.}$ equivariant Borel map defined in \eqref{eq:xi*b}. Then $\xi _\mathsf{b} ^*$ is countable-to-one $\nu \textsl{-a.e.}$, and the $\alpha $-invariant probability measure $\mu =(\xi _\mathsf{b} ^*)_*\nu $ on $X$ is singular with respect to Haar measure and satisfies that $h_\nu (\bar{\sigma })=h_{\mu }(\alpha )$.
\end{theo}

For the proof of Theorem \ref{t:weaklybounded} we need several lemmas. The hypotheses of these lemmas are those of the theorem.

We call two points $v,v'\in V_\beta $ \emph{equivalent} (in symbols: $v \sim v'$) if $v-v'\in f(\bar{\sigma })(\ell ^ *(\mathbb{Z},\mathbb{Z}))$ or, equivalently, if $\xi ^*(v)-\xi ^*(v')\in X ^{(0)}$ (cf. \eqref{eq:inverse2}). Denote by
\begin{equation}
\label{eq:R}
\mathbf{R}=\{(v,v'):v \sim v'\} \subset V_\beta \times V_\beta
\end{equation}
the resulting equivalence relation, and write
\begin{equation}
\label{eq:class}
\mathbf{R}(v)=\{ v'\in V_\beta :v \sim v'\}
\end{equation}
for the equivalence class of every $v \in V_\beta $.

\begin{lemm}
\label{l:Borel}
The set $\mathbf{R}\subset V_\beta \times V_\beta $ is Borel and $\bar{\sigma }\times \bar{\sigma }$-invariant.
\end{lemm}

\begin{proof}
For every $r>0$, the sets
$$
B_r ^*=\biggl\{ v \in \ell ^*(\mathbb{Z},\mathbb{Z}):\sup_{n \in \mathbb{Z}} \,\frac{|v_n|}{|n|+1}\le r\biggr\}
$$
and $C_r=f(\bar{\sigma })(B_r ^*)\subset \ell ^*(\mathbb{Z},\mathbb{Z})$ are compact, and the map $p\colon C_r \times V_\beta \longrightarrow \ell ^ *(\mathbb{Z},\mathbb{Z})\times V_\beta $, given by $p(v',v)=(v'+v,v)$, is continuous. Hence $\tilde{C}_r=p(C_r \times V_\beta )\cap (V_\beta \times V_\beta )$ is compact and $\mathbf{R}=\bigcup_{r>0}\tilde{C}_r$ is Borel. The $\bar{\sigma }\times \bar{\sigma }$-invariance of $\mathbf{R}$ is obvious.
\end{proof}

For every subset $F \subset \mathbb{Z}$ we write $\pi _F\colon \ell ^ \infty (\mathbb{Z},\mathbb{Z})\longrightarrow \mathbb{Z}^ F$ for the projection onto the coordinates in $F$.

%%%%CHANGES
\begin{lemm}
\label{l:finite}
Let $Y \subset V_\beta $ be a shift-invariant Borel set with $\nu (Y)=1$ such that \eqref{eq:b1} holds for every $v \in Y$, and let
\begin{equation}  \label{equation with three definitions}
\begin{aligned}
Y(M)& =\{ y \in Y:\| \mathsf{b}(y)\|_\infty \le M \},
\\
L(M)& =\{ y \in \ell ^*(\mathbb{Z},\mathbb{Z}):\| y-\bar{\xi }^*\circ f(\bar{\sigma })(y)\|_\infty \le M \}
\\
\mathbf{R}(M,w)& =\bigl(w+f(\bar{\sigma })(L(M))\bigr)\cap V_\beta \subset \mathbf{R}(w)
\end{aligned}
\end{equation}
for every $M\ge1$ and $w \in V_\beta $ \textup{(}cf. \eqref{eq:inverse2}\textup{)}. Then for every $w \in V_\beta $
\begin{equation}
\label{eq:trivial}
R(w) = \bigcup_ {M=1}^ \infty R(M,w)
\end{equation}
and for every $K,M,n\ge1$
\begin{equation}
\label{eq:finite}
\begin{aligned}
\absolute {\pi _{\{ 0,\dots ,n \}}\bigl(\mathbf{R}(K,w) \cap Y(M)\cap \bar \sigma ^{-n}(Y(M))\bigr)} &\le c(M,K),
\\
\absolute {\pi _{\{-n,\dots ,0 \}}\bigl(\mathbf{R}(K,w) \cap Y(M)\cap \bar{\sigma }^ n(Y(M))\bigr)}
&\le c(M,K),
\end{aligned}
\end{equation}
where $c(M,K)$ is a constant depending only on $K$, $M$ and $f$.
\end{lemm}

\begin{proof}
We first prove \eqref{eq:trivial}. Indeed, by \eqref{eq:inverse2} and the remarks following $ \bigcup_ M L(M) = \ell ^{*} (\ZZ, \ZZ)$, hence
\begin{equation*}
\bigcup_ {M = 1 } ^ \infty R (M, w) = \left (w + f(\bar \sigma) (\ell ^{*} (\ZZ, \ZZ)) \right) \cap V _ \beta = R (w)
.\end{equation*}

We now turn to prove \eqref{eq:finite}.
By \eqref{eq:xi*} there exists a constant $M_1>0$ such that
$$
\max_{j=0,\dots ,m}|\bar{\xi }^*(w)_j|\le M_1 \cdot \| w \|_\infty
$$
for every $w \in \ell ^ \infty (\mathbb{Z},\mathbb{Z})$. As $(\bar{\sigma }^*)^ n \circ \bar{\xi }^*(w)=\bar{\xi }^*\circ \bar{\sigma }^ n(w)+\mathsf{b}(\bar{\sigma }^ nw)-\bar{\sigma }^ n \mathsf{b}(w)$ for every $n \in \mathbb{Z}$	,
$$
\max_{j=0,\dots ,m}|\bar{\xi }^*(w)_{n+j}|\le M_1 \beta +2M
$$
for every $M\ge1$, $n \in \mathbb{Z}$ and $w \in Y(M)\cap \bar{\sigma }^{-n}(Y(M))$. We fix $w \in V_\beta $ and obtain that, for every $v \in \mathbf{R}(K,w)\cap Y(M)\cap \bar{\sigma }^{-n}(Y(M))$
$$
\max_{j=0,\dots ,m}|\bar{\xi }^*(v)_j|\le M_1 \beta ,\enspace \enspace \max_{j=0,\dots ,m}|\bar{\xi }^*(v)_{n+j}|\le M_1 \beta +2M,
$$
and that there exists a unique $y \in \ell ^*(\mathbb{Z},\mathbb{Z})$ with $v=w+f(\bar{\sigma }^*)(y)$ and $\| y-\bar{\xi }^*\circ f(\bar{\sigma }^*)(y)\|_\infty \le K$.

If $v'$ is a second element in $\mathbf{R}(K,w)\cap Y(M)\cap \bar{\sigma }^{-n}(Y(M))$ with $v'=w+f(\bar{\sigma }^*)(y')$ for some $y'\in \ell ^*(\mathbb{Z},\mathbb{Z})$, then $\| y'-\bar{\xi }^*\circ f(\bar{\sigma }^*)(y')\|_\infty \le K$, and hence
$$
\max_{j=0,\dots ,m}|y_j-y_j'|\le 2M_1 \beta +2K\enspace \textup{and}\enspace \max_{j=0,\dots ,m}|y_{n+j}-y_{n+j}'|\le 2M_1 \beta +4M+2K.
$$
For every $L>0$ we set
\begin{align*}
B(n,L)&=\{ w \in \ell ^* (\mathbb{Z},\mathbb{Z}):|w|_j\le L
\\
&\enspace \enspace \enspace \textup{for}\; 0\le j\le m\enspace \textup{and}\enspace n\le j\le n+m \},
\end{align*}

If the first inequality in \eqref{eq:finite} does not hold for some $w \in V_\beta $ and $n>0$, then we can find elements $y,z \in B(n,2M_1 \beta +4M+2K+1)$ with the following properties:
\begin{gather*}
(y_0,\dots ,y_m)=(z_0,\dots ,z_m),
\\
(y_n,\dots ,y_{n+m})=(z_n,\dots ,z_{n+m}),
\\
(y_{m+1},\dots ,y_{n-1})\ne (z_{m+1},\dots ,z_{n-1}).
\end{gather*}
so that $\bar{y}=w+f(\bar{\sigma }^*)(y)$, $\bar{z}=w+f(\bar{\sigma }^*)(z)$ are both in $V_\beta$. Note that these two points satisfy $\bar{y}_0=\bar{z}_0$, $\bar{y}_n=\bar{z}_n$ and $(\bar{y}_1,\dots ,\bar{y}_{n-1})\ne (\bar{z}_1,\dots ,\bar{z}_{n-1})$.

Suppose, without loss in generality, that
$$
(\bar{y}_1,\dots ,\bar{y}_{n-1})\prec (\bar{z}_1,\dots ,\bar{z}_{n-1})
$$
and hence
$$
(\bar{y}_1,\bar{y}_2,\dots )\prec (\bar{z}_1,\bar{z}_2,\dots )
$$
in the lexicographic order. We set
$$
y'_j=
\begin{cases}
z_j&\textup{if}\enspace j\le 0,
\\
y_j&\textup{if}\enspace j>0,
\end{cases}
$$
and put $z'=z$. Then $\bar{y}'=w+f(\bar{\sigma })(y')$ is of the form
$$
\bar{y}_j'=
\begin{cases}
\bar{z}_j&\textup{if}\enspace j\le 0,
\\
\bar{y}_j&\textup{if}\enspace j>0,
\end{cases}
$$
and $\bar{y}'\in V_\beta $ by Proposition \ref{p:beta} (1). Put $\bar{z}'=w+f(\bar{\sigma }^*)(z')=\bar{z}$, remember that $\bar{y}_j'=\bar{z}_j'$ for $j\le 0$ and for $j=n$, and assume for the moment that
$$
(\bar{y}'_{n+1},\bar{y}'_{n+2},\dots )\prec (\bar{z}'_{n+1},\bar{z}'_{n+2},\dots )
$$
in the lexicographic order (if this is not the case we have to interchange the roles of $y'$ and $z'$ below). Let
$$
z_j''=
\begin{cases}
z'_j&\textup{if}\enspace j\le n,
\\
y'_j&\textup{if}\enspace j>n,
\end{cases}
$$
and set $y''=y'$ and $\bar{y}''=\bar{y}'$. The point $\bar{z}''=w+f(\bar{\sigma }^*)(z'')$ is of the form
$$
\bar{z}_j''=
\begin{cases}
\bar{z}'_j&\textup{if}\enspace j\le n,
\\
\bar{y}'_j&\textup{if}\enspace j>n,
\end{cases}
$$
and lies in $V_\beta $ by Proposition \ref{p:beta} (1). By construction, $y_j''=z_j''$ for $j\le m$ and $j\ge n$, and hence $v''=y''-z''\in \ell ^ 1(\mathbb{Z},\mathbb{Z})$. Since $\bar{z}''$ and $\bar{y}''=\bar{z}''+f(\bar{\sigma })v''$ lie in $V_\beta $ we obtain a contradiction to Proposition \ref{p:beta} (2). This proves the first inequality in \eqref{eq:finite}, and the proof of the second one is analogous.
\end{proof}

\begin{lemm}
\label{l: finite 2}
Let $R(K,w), Y(M), c(M,K)$ be as in Lemma~\ref{l:finite}, and let
\begin{equation*}
\tilde Y (M) = Y (M) \cap \left\{ y \in Y: \liminf_ {n \to \pm \infty} \frac {1 }{ \absolute n} \sum_ {k = n } ^ {2 n} 1 _ {Y (M)} (\bar \sigma ^ k w) \geq \frac {1 }{ 2} \right\}
.\end{equation*}
Then for every $K, M$ and $w \in V _ \beta$
\begin{equation} \label{equation needed to avoid summations}
\absolute {R(K, w) \cap \tilde Y (M)} \leq 100c(M,2K) ^2
\end{equation}
\end{lemm}

\begin{proof}
%We first note that by the definitions in
%(equation: equation with three definitions) for any "w in V_Greek beta" and any "w' in R (K, w)"
%(begin an equation)
%R (K, w) subset R (2K, w')
%(end this equation).
Assume in contradiction that there is some $w \in V _ \beta$ for which \eqref{equation needed to avoid summations} fails. Then there is a $n_0$ so that at least one of the following holds:
\begin{align*}
\absolute {\pi _ {\left\{ 0, n \right\}} \left (R(K, w) \cap \tilde Y (M) \right)} & > 10c(M,2K) \qquad \text{for every $n>n_0$ or} \\
\absolute {\pi _ {\left\{ - n, 0 \right\}} \left (R(K,w) \cap \tilde Y (M ) \right)} & > 10c(M,2K) \qquad \text{for every $n>n_0$}
.\end{align*}
Assume that the former holds (the argument for the latter is identical). Suppose $w_1$, \dots, $w_{10c(M,2K)}$ are $10c(M,2K)$ points in $R(K,w) \cap \tilde Y (M)$ with $ \pi _ {\left\{ 0, n \right\}} (w _ i) \neq \pi _ {\left\{ 0, n \right\}} (w _ j) $
for $i \ne j$. Then by definition of $\tilde Y (M)$, for $n_1 > n_0$ sufficiently large
\begin{equation*}
\sum_ {i = 1} ^ {10 c (M, 2 K) + 1} \sum_ {k = n _ 1} ^ {2 n _ 1} 1 _ {Y (M)} (\bar \sigma ^ k w) > 4 n_1 c (M, 2 K)
\end{equation*}
so that there would be some $n_2 > n_0$ for which at least $c (M, 2 K) +1$ of the $w_i$ (which without loss of generality we can assume to be $w _ 1, \dots, w _ {c (M, 2 K) +1} $)
satisfy $\bar \sigma ^ {n_2} w_i \in Y(M)$.
We already know all the $w_i$ are in $R(K,w) \cap \tilde Y(M) \subset R(K,w) \cap Y(M)$. Since $\pi _ {\left\{ 0, n_0 \right\}} w_i$ are all distinct (which also implies that $\pi _ {\left\{ 0, n_2 \right\}} w_i$ are all distinct), the points $w_1$, \dots $w_{c (M, 2 K) +1}$ show that
\begin{equation*}
\absolute {\pi _ {\left\{ 0,n_2 \right\}} \left (R(K,w) \cap Y (M) \cap \bar \sigma ^ {-n_2} Y(M) \right)} \geq c (M, 2 K) +1
\end{equation*} which is in contradiction to \eqref{eq:finite}.
\end{proof}

\begin{proof}[Proof of Theorem \ref{t:weaklybounded}]

Let $\nu$ and $\mu$ be measures on $V _ \beta$ and $X$ respectively as in Theorem~\ref{t:weaklybounded}. We will show in fact something stronger than merely that $\xi _\mathsf{b} ^*$ is countable-to-one: we will show that there is a subset $Z_1 \subset V _ \beta$%%%CHANGE
with $\nu (Z_1)=1$ so that for any $x \in X$,
\begin{equation*}
[\xi _ \mathsf{b} ^ {*}] ^{-1} \left (x +X^{(0)} \right) \cap Z_1
\end{equation*}
is countable.
Indeed, take $Z_1 = \bigcup_ {M = 1} ^ \infty \tilde Y (M)$, with $\tilde Y (M)$ as in Lemma~\ref{l: finite 2}; clearly $\nu (Z_1) = 1$. For any $x = \xi _ \mathsf{b} ^ {*} (w) \in \xi _ \mathsf{b} ^ {*} (V_\beta)$
\begin{equation*}
[\xi _ \mathsf{b} ^ {*}] ^{-1} \left (x +X^{(0)} \right) \cap Z_1 = R(w) \cap Z_1 = \bigcup_ {K, M = 1 } ^ \infty (R(K,w) \cap \tilde Y (M))
.\end{equation*}
By Lemma~\ref{l: finite 2}, $R(K,w) \cap \tilde Y (M)$ is finite and the result follows.

Since countable-to-one factor maps do not decrease entropy, $h_\nu (\bar{\sigma })=h_{\mu }(\alpha )$. Furthermore, the set $Z_2=\xi _\mathsf{b} ^*(Z_1)\subset X$ satisfies $\mu (Z_2)=1$ and intersects each coset of $X ^{(0)}$ in a countable set. Hence by Fubini $\lambda _X(Z_2)=0$, which proves that $\lambda _X$ and $\mu $ are mutually singular.
\end{proof}

As we have seen, on $V _ \beta$ there is a unique $\bar \sigma$-invariant measure $\mu _ \beta$ with maximal entropy $\log \beta$. If this measure would have been weakly $\mathsf{d}$-bounded, $[\xi _ \mathsf{b} ^ {*}]_* \mu _ \beta$ would have been a measure on $X$ which has entropy $\log \beta$ but is singular with respect to $\lambda _ X$, which is clearly absurd as $\lambda _ X$ is the unique $\alpha$-invariant measure on $X$ with entropy $\log \beta$. Thus as a biproduct of our discussion on symbolic representations we obtain:

\begin{coro}
The measure $\mu _ \beta$ on $V _ \beta$ is not weakly $\mathsf{d}$-bounded.
\end{coro}

\section {Some examples of invariant measures in the Salem case}\label{s: examples}

\begin{theo}
\label{t:salem2}
Let $\beta >1$ be a Salem number, and let $V_\beta \subset \ell ^ \infty (\mathbb{Z},\mathbb{Z})$ be the corresponding two-sided beta-shift space. For every $\varepsilon >0$ there exists a $\mathsf{d}$-bounded shift-invariant probability measure $\nu $ on $V_\beta $ with $h_\nu (\sigma _\beta )>\log \beta -\varepsilon $, where $\sigma _\beta =\bar{\sigma }_{V_\beta }$ is the beta-shift.
\end{theo}

\begin{proof}
As in the proof of Theorem \ref{t:3} we choose an enumeration $\Omega _f ^{(0)} =\{ \omega _1,\dots ,\omega _{m_0}\}$ of $\Omega _f ^{(0)} $, write $\mathbf{W}_f ^{(0)} = \mathbb{C}\otimes_\mathbb{R}W_f ^{(0)}$ for the complexification of $W_f ^{(0)}$, and use the basis $\{ w(\omega _i):i=1,\dots ,{m_0}\}$ in \eqref{eq:womega} to identify $\mathbf{W}_f ^{(0)}$ with $\mathbb{C}^{m_0}$.
% EBL: fixed minor error
Let
\begin{equation*}
\Gamma _ \beta = \overline {\left\{ (\omega _ 1 ^ n, \dots, \omega _ {m _ 0} ^ n): n \in \ZZ \right\}}
.
\end{equation*}
Put $\mathbf{V}=V_\beta \times \Gamma _ \beta$ (cf. \eqref{eq:Sm0}) and define a map $S_\beta \colon \mathbf{V}\longrightarrow \mathbf{V}$ by $S_\beta (v,\gamma )=(\bar{\sigma }v,M_{\boldsymbol{\omega }}\gamma )$ for every $v \in V_\beta $ and $\gamma =(\gamma _1,\dots ,\gamma _{m_0})\in \Gamma _ \beta \subset \mathbb{C}^{m_0}$, where $M_{\boldsymbol{\omega }}$ is defined in \eqref{eq:Mgamma} and \eqref{eq:pmb}.

Let $\lambda $ be the Haar (= normalized Lebesgue) measure on $\Gamma _ \beta$. Since the unique shift-invariant measure of maximal entropy $\mu _\beta $ on $V_\beta $ is mixing (cf. \cite{Hofbauer}), the measure $\mu _\beta \times \lambda $ on $\mathbf{V}$ is ergodic under $T_\beta $. The map $F_\beta \colon \mathbf{V}\longrightarrow \mathbb{C}^{m_0}$, given by $F_\beta (v,\gamma )_i=\gamma _iv_0$ for every $v=(v_n)\in V_\beta $, $\gamma =(\gamma _1,\dots ,\gamma _{m_0})\in \Gamma _ \beta$ and $i=1,\dots ,m_0$, satisfies that $\int F_\beta \,d(\mu _\beta \times \lambda )=0$, and the ergodic theorem implies that
$$
\lim_{K \to \infty }\biggl\| \frac 1K \sum_{k=0}^{K-1}F_\beta (T_\beta ^ k(v,\gamma ))\biggr\|_\infty =\lim_{K \to \infty }\biggl|\frac 1K \sum_{k=0}^{K-1}\omega _i ^ kv_k\biggr|=0\enspace ( \mu _\beta \times \lambda )\textsl{-a.e.}
$$
for $i=1,\dots ,m_0$, where $\| \cdot \|_\infty $ is the maximum norm on $\mathbb{C}^{m_0}$. Hence

$$
\lim_{K \to \infty }\| \mathsf{d}(K,v)/K \|_\infty =0\enspace \mu _\beta \textsl{-a.e.}
$$
We fix a positive integer $J$ and choose $K>0$ sufficiently large so that $\mu _\beta (B_{K,J})>1-1/J$, where
$$
B_{K,J}=\{ v \in V_\beta :\| \mathsf{d}(k,v)\|_\infty \le K\enspace \textup{for}\enspace k=0,\dots ,KJ \}.
$$
Note that the set $B_{K,J}$ is a union of cylinder sets which depend only on the coordinates $0,\dots ,KJ-1$.

Since $M_{\boldsymbol{\omega }}$ acts minimally on $\Gamma _ \beta$, there exists an $L>0$ with the following property: for every pair $\mathbf{v}, \mathbf{w}\in \mathbb{C}^{m_0}$ of vectors with $\| \mathbf{v}\|_\infty \le 1$ and $\| \mathbf{w}\|_\infty \le 1$ there exists an $l \in \{ 0,\dots ,L-1 \}$ with $\| M_{\boldsymbol{\omega }}^ l \mathbf{v} +\mathbf{w}\|_\infty \le 1$.

Let $v \in V_\beta $. By inserting zero coordinates in an appropriate manner we modify $v$ to a point $v ^*\in V_\beta $ with $v_n ^*=v_n$ for $n<0$ such that
% EBL changed epsilon to four
$\| \mathsf{d}(m,v ^*)\|_\infty <4 K$ for every $m\ge0$.

In order to describe this modification we proceed by induction and assume that $v=v ^{(0)}\in V_\beta $. If $v \in B_{K,J}$ we put $v ^{(1)}=v$ and $\mathbf{v}^{(1)}=\mathsf{d}(JK,v ^{(0)})$.

If $v\notin B_{J,K}$ we use our choice of $L$ to find an integer $l_1 \in \{ 0,\dots ,L-1 \}$ such that the point $v(1)$, given by
$$
v(1)_n=
\begin{cases}
v_n&\textup{if}\enspace n\le K-1,
\\
0&\textup{if}\enspace n=K,\dots ,K+l_1-1,
\\
v_{n-l_1}&\textup{if}\enspace n\ge K+l_1,
\end{cases}
$$
which satisfies that $\| \mathsf{d}(2K+l_1,v(1))\|_\infty \le K$. Next we choose $l_2 \in \{ 0,\dots ,L-1 \}$ such that the point $v(2)$ with
$$
v(2)_n=
\begin{cases}
v(1)_n&\textup{if}\enspace n\le 2K+l_1-1,
\\
0&\textup{if}\enspace n=2K+l_1,\dots ,2K+l_1+l_2-1,
\\
v_{n-l_1-l_2}&\textup{if}\enspace n\ge 2K+l_1+l_2,
\end{cases}
$$
satisfies that $\| \mathsf{d}(3K+l_1+l_2,v(2))\|_\infty \le K$. By continuing in this manner we eventually obtain integers $l_1,\dots ,l_{J-1}\in \{ 0,\dots ,L-1 \}$ and a point $v ^{(1)}=v(J-1)\in V_\beta $ (cf. Proposition \ref{p:beta} (1)) with
$$
v ^{(1)}_n=
\begin{cases}
v_n&\textup{if}\enspace n\le K-1,

\\
0&\textup{if}\enspace n=K,\dots ,K+l_1-1,
\\
\vdots
\\
0&\textup{if}\enspace n=(J-1)K+l_1+\dots +l_{J-2}, \dots ,
\\
&\qquad \qquad(J-1)K+l_1+\dots +l_{J-1}-1,
\\
v_{n-l_1-\dots -l_{J-1}}&\textup{if}\enspace n\ge (J-1)K+l_1+\dots +l_{J-1},
\end{cases}
$$
satisfies that $\| \mathsf{d}(JK+l_1+\dots +l_{J-1},v ^{(1)})\|_\infty \le K$. We set $l ^{(1)}=l_1+\dots +l_{J-1}$, $\mathbf{v}^{(1)}=\mathsf{d}(JK+l ^{(1)},v ^{(1)})$ and note that $\bar{\sigma }^{JK+l ^{(1)}}v ^{(1)}\in B_{K,J}$ if and only if $\bar{\sigma }^{JK}v \in B_{K,J}$, and that
$$
\| \mathsf{d}(j,v ^{(1)})\|_\infty \le 2K
$$
for $j=0,\dots ,JK+l ^{(1)}$.

We repeat this process with $v$ replaced by $w=\bar{\sigma }^{JK+l ^{(1)}}v ^{(1)}$ and obtain an integer $l ^{(2)}\in \{ 0,\dots ,J(L-1)\}$ and a point $w'\in V_\beta $ with the following properties.
\begin{enumerate}
\item[(i)]
$\| \mathsf{d}(JK+l ^{(2)},w')\|_\infty \le K$ and $\| \mathsf{d}(j,w')\|_\infty \le 2K$ for $j=0,\dots ,JK+l ^{(2)}$,
\item[(ii)]
$w'_n=w_n$ for $n<0$ and $w'_{n+l ^{(2)}}=w_n$ for $n\ge JK$,
\item[(iii)]
$w'$ is obtained from $w$ by inserting $l ^{(2)}\le (J-1)(L-1)$ zeros among the coordinates $w_0,\dots ,w_{JK-1}$, and $l ^{(2)}=0$ if and only if $\bar{\sigma }^{JK}(v)\in B_{K,J}$.
\end{enumerate}
Next we set $\mathbf{w}=\mathsf{d}(JK+l ^{(2)},w')$, choose a $j ^{(1)}\in \{ 0,\dots ,L-1 \}$ with $\| M_{\boldsymbol{\omega }}^{j ^{(1)}}\mathbf{v}^{(1)}+\mathbf{w}\|_\infty \le K$, and define $v ^{(2)}\in V_\beta $ by
$$
v ^{(2)}_n=
\begin{cases}
v ^{(1)}_n&\textup{if}\enspace n<JK+l ^{(1)},
\\
0&\textup{if}\enspace n=JK+l ^{(1)},\dots ,JK+l ^{(1)}+j ^{(1)}-1
\\
w'_{n-JK-l ^{(1)}-j ^{(1)}}&\textup{if}\enspace n\ge JK+l ^{(1)}+j ^{(1)}.
\end{cases}
$$
The point $v ^{(2)}$ lies in $V_\beta $ by Proposition \ref{p:beta} (1) and has the following properties.
\begin{enumerate}
\item[(i')]
$\| \mathsf{d}(2JK+l ^{(1)}+j ^{(1)}+l ^{(2)}, v ^{(2)})\|_\infty \le K$ and $\| \mathsf{d}(j,v ^{(2)})\|_\infty \le 2K$ for $j=0,\dots ,2JK+l ^{(1)}+j ^{(1)}+l ^{(2)}$,
\item[(ii')]
$v ^{(2)}_n=v_n$ for $n<0$ and $v ^{(2)}_{n+l ^{(1)}+j ^{(1)}+l ^{(2)}}=v_n$ for $n\ge 2JK$,
\item[(iii')]
$v ^{(2)}$ is obtained from $v ^{(1)}$ by inserting $l ^{(2)}\le J(L-1)$ zeros among the coordinates $v_{JK+l ^{(1)}},\dots ,v_{2JK-1+l ^{(1)}}$, and $l ^{(2)}=0$ if and only if $\bar{\sigma }^{JK}v \in B_{K,J}$ (or, equivalently, if and only if $\bar{\sigma }^{JK+l ^{(1)}}v ^{(1)}\in B_{K,J}$).
\end{enumerate}

By repeating this process we obtain sequences $(v ^{(m)},\,m\ge1)$ in $V_\beta $ and $(l ^{(m)},\,m\ge1)$ and $(j ^{(m)},\,m\ge1)$ of positive integers satisfying the following conditions for every $m\ge1$.
\begin{enumerate}
\item
$0\le l ^{(m)}\le J(L-1)$ and $0\le j ^{(m)}\le L-1$,
\item
If $L ^{(m)}=\sum_{i=1}^ m l ^{(i)}$, $J ^{(m)}=\sum_{i=1}^{m-1}j ^{(i)}$ and $L ^{(i)}=J ^{(i)}=0$ for $i\le0$, then
$$
\smash[t]{L ^{(m)}\le J(L-1)\cdot \sum_{i=0}^{m-1}1_{V_\beta \smallsetminus B_{K,J}}(\bar{\sigma }^{iJK}v),}
$$
where $1_S$ denotes the indicator function os a set $S \subset V_\beta $, and
$$
\| \mathsf{d}(mJK+L ^{(m)}+J ^{(m-1)},v ^{(m)})\|_\infty \le K,
$$
\item
$\| \mathsf{d}(j,v ^{(m)})\|_\infty \le 2K\enspace \textup{for}\enspace j=0,\dots ,mJK+L ^{(m)}+J ^{(m-1)}$, \vspace{4mm}
\item
~\vspace{-9.8mm}
$$
\qquad v ^{(m)}_n=
\begin{cases}
v_n ^{(m-1)}&\textup{if}\enspace n<(m-1)JK+L ^{(m-1)}+J ^{(m-2)},
\\
v_{n-L ^{(m)}-J ^{(m-1)}}&\textup{if}\enspace n\ge mJK+L ^{(m)}+J ^{(m)}.
\end{cases}
$$
\end{enumerate}
From the conditions (3)--(4) above we see that the sequence $(v ^{(m)},\,m\ge1)$ converges to an element $v ^*\in V_\beta $ with
$$
\| \mathsf{d}(j,\bar{\sigma }^{j'}v ^*)\|_\infty \le 4K
$$
for every $j,j'\ge0$.

If $m$ is sufficiently large, then the set
$$
C_m=\biggl\{ v \in V_\beta :\sum_{i=0}^{m-1}1_{V_\beta \smallsetminus B_{K,J}}(\bar{\sigma }^{iJK}v)\le 2m/J\biggr\}
$$
has $\mu _\beta $-measure $>1-1/J$.

So far we have kept $J$ and $K$ fixed, but now we begin to vary them. If
$$
P(m)=\pi _{\{ 0,\dots ,mJK-1 \}}(C_m)
$$
is the projection of the set $C_m$ onto the coordinates $0,\dots ,mJK-1$, then the Shannon-McMillan-Breiman theorem applied to $\mu _ \beta$ (cf. \cite{Parrybook}) implies that the cardinality of $P(m)$ satisfies that
$$
\lim_{J \to \infty }\frac 1{mJK}\log\,|P(m)|=\log\,\beta ,
$$
since $h_{\mu _\beta }(\bar{\sigma }_\beta )=\log\,\beta $ (note that $K$ depends on $J$ and tends to infinity as $J \to \infty $). We fix $\varepsilon >0$ and choose $J$ (and hence $K$) sufficiently large so that $P(m)>(\beta -\varepsilon )^{mJK}$ for all sufficiently large $m$. For every $v \in C_m$, the number of zero coordinates inserted among the coordinates $v_0,\dots ,v_{mJK-1}$ in the transition from $v$ to $v ^*$ is less than $m \cdot (L-1)+2m \cdot (L-1)\cdot K/J$, so that
$$
|\pi _{\{ 0,\dots ,mJK-1 \}}(\{ v ^*:v \in C_m \})|\ge |\pi _{\{ 0,\dots ,m \cdot (JK-L-2 \cdot (L-1)\cdot K/J)\}}(C_m)|.
$$
This shows that, for sufficiently large $K$, the topological entropy of the closed, $\bar{\sigma }$-invariant subset
$$
\{ v \in V_\beta :\| \mathsf{d}(j,\bar{\sigma }^{j'}v ^*)\|_\infty \le 4K\enspace \textup{for every}\enspace j\ge0\enspace \textup{and}\enspace j'\in \mathbb{Z}\}
$$
is arbitrarily close to $\log\,\beta $, and the variational principle (cf. \cite{Walters}) guarantees that we can find $\bar{\sigma }$-invariant and ergodic probability measures $\nu $ on $V_\beta $ with entropy arbitrarily close to $\log\,\beta $.
\end{proof}

\section*{Acknowledgment}
This research has been supported in part by NSF grant DMS 0140497
(E.L.) and FWF Project P16004--N05 (K.S.). During part of this work,
both authors received support from the American Institute of
Mathematics and NSF grant DMS 0222452.  We would furthermore like to
express our gratitude to the Mathematics Departments of the
University of Washington, Stanford University, the Newton
Institute, Cambridge and the ETH Z\"urich for hospitality during
parts of this work. E.L. is a Clay Research Fellow and is grateful for this generous support from the Clay Mathematics Institute. E.L. would also like to thank Rick Kenyon for an interesting and helpful discussion on these and related topics.

\end{document}